\documentclass[11pt,a4paper]{article}
\usepackage[T1]{fontenc}
\usepackage[utf8]{inputenc}
\usepackage{lmodern}
\usepackage{amsmath,amssymb,amsthm,mathtools}
\usepackage[margin=28mm]{geometry}
\usepackage{microtype}
\usepackage{enumitem}
\usepackage[hidelinks]{hyperref}
\usepackage{bookmark}
\allowdisplaybreaks[2]
\numberwithin{equation}{section}
\newtheorem{theorem}{Theorem}
\newtheorem{lemma}{Lemma}[section]
\newtheorem{corollary}{Corollary}[section]
\theoremstyle{definition}
\newtheorem{definition}{Definition}[section]
\newtheorem{example}{Example}[section]
\theoremstyle{remark}
\newtheorem{remark}{Remark}[section]
\newcommand{\R}{\mathbb R}
\newcommand{\N}{\mathbb N}

\newcommand{\Gat}{G\^ateaux}
\newcommand{\Fre}{Fr\'echet}
\newcommand{\supp}{\operatorname{supp}}
\newcommand{\Ls}{\mathcal L_{\mathrm s}}
\newcommand{\norm}[1]{\left\lVert #1\right\rVert}
\newcommand{\tnorm}[1]{\left\lvert\left\lvert\left\lvert #1\right\rvert\right\rvert\right\rvert}
\newcommand{\abs}[1]{\left|#1\right|}
\newcommand{\dd}{\,\mathrm d}
\hypersetup{
  pdftitle={No LUR Renorming of ℓp Is Three Times Gâteaux Differentiable for 2 < p < 3},
  pdfauthor={Diana Nedelcheva and Boyan Zlatanov},
  pdfsubject={Generic local regularity, block estimates, and a geometric obstruction for 2<p<3},
  pdfkeywords={Gâteaux differentiability, locally uniformly rotund norm, Kadec norm, strong exposure, bump function, block sequence}
}
\title{No LUR Renorming of \(\ell_p\) Is Three Times\\
\Gat\ Differentiable for \(2<p<3\)}
\author{Diana Nedelcheva\thanks{Department of Mathematics and Physics,
  Technical University of Varna, 1 Studentska Str., 9000 Varna, Bulgaria.
  \newline \href{mailto:diana.nedelcheva@tu-varna.bg}{\nolinkurl{diana.nedelcheva@tu-varna.bg}}
  \enspace\textperiodcentered\enspace
  \href{https://orcid.org/0000-0001-6094-7595}{ORCID 0000-0001-6094-7595}}
  \and
  Boyan Zlatanov\thanks{Department of Mathematical Analysis, Faculty of
  Mathematics and Informatics, University of Plovdiv Paisii Hilendarski,
  4000 Plovdiv, Bulgaria.
  \newline Correspondence: \href{mailto:bobbyz@uni-plovdiv.bg}{\nolinkurl{bobbyz@uni-plovdiv.bg}}
  \enspace\textperiodcentered\enspace
  \href{https://orcid.org/0000-0002-5857-9372}{ORCID 0000-0002-5857-9372}}}
\date{November 29, 2026}

\begin{document}
\begingroup
\makeatletter
\let\@fnsymbol\@arabic 
\makeatother
\maketitle
\endgroup
\begin{abstract}
For \(2<p<3\), no equivalent locally uniformly rotund (LUR) norm on
\(\ell_p\) is three times \Gat\ differentiable at every nonzero point
in the recursive operator-norm sense. We prove this for the larger
class of Kadec norms. For any equivalent norm with this differentiability,
the strongly exposed points of its unit ball are nowhere dense in its
unit sphere. The proof uses the fact that a continuous real-valued
function with three such derivatives on an open subset of a real Banach
space has a locally Lipschitz second \Fre\ derivative on a dense open
subset of its domain. We also construct a convex function \(\Psi\) on
\(\ell_p\), twice \Fre\ differentiable everywhere, whose third \Gat\
derivative vanishes at zero. Its quadratic Taylor polynomial \(P_2\) at
zero satisfies \(\sup_n|\Psi(tu_n)-P_2(tu_n)|=|t|^p\) for every normalized
successive block sequence \((u_n)\) and every \(t\ne0\).
\end{abstract}



\section{Introduction}\label{sec:intro}

Let \(2<p<3\). We prove that no equivalent Kadec norm on \(\ell_p\) is
three times \Gat\ differentiable at every nonzero point in the recursive
operator-norm sense (Theorem~\ref{thm:kadec}). In particular, no equivalent
locally uniformly rotund (LUR) norm has this differentiability. More precisely, for
any equivalent norm with this differentiability, the strongly exposed
points of its unit ball are nowhere dense in its unit sphere. Our argument
leaves open whether \(\ell_p\) admits an equivalent norm with this
differentiability without the Kadec assumption (Remark~\ref{rem:remaining}).

Our second result concerns the uniformity of Taylor remainders at a
single point (Theorem~\ref{thm:counter}). We construct a convex function
\(\Psi\) on \(\ell_p\) with \(D^3\Psi(0)=0\). Along each fixed direction
the remainder after its quadratic Taylor polynomial is \(o(|t|^3)\).
Yet for every normalized successive block sequence \((u_n)\), the supremum
of the absolute remainders at \(tu_n\) is exactly \(|t|^p\).
Thus even passing to a block subsequence cannot give a common cubic bound.

We use the distinction between the smoothness of a given norm and the
geometry of its unit ball throughout
\cite{Day,Cudia,DeBernardiPretiSomaglia}. We first specify our norm and
derivative conventions, since several notions of higher-order \Gat\
differentiability occur in the literature
\cite{MaleevZlatanov,McLaughlinVanderwerff,ITZ}.

\subsection{The original norm, an equivalent norm, and their dual norms}
\label{subsec:norm-conventions}

Throughout the paper, \(X\) is a real Banach space. Its original norm
is denoted by \(\norm{\cdot}\), and an equivalent norm is denoted by
\(\tnorm{\cdot}\). Thus there are constants \(a,b>0\) such that
\begin{equation}\label{eq:norm-equivalence-intro}
 a\norm{x}\le\tnorm{x}\le b\norm{x}\qquad(x\in X).
\end{equation}
Unless explicitly stated otherwise, \(\norm{\cdot}^{*}\) and
\(\tnorm{\cdot}^{*}\) denote the norms on the continuous dual space
\(X^*\) induced by \(\norm{\cdot}\) and \(\tnorm{\cdot}\),
respectively:
\begin{equation}\label{eq:dual-norms-intro}
 \norm{f}^{*}=\sup_{\norm{x}\le1}|f(x)|,\qquad
 \tnorm{f}^{*}=\sup_{\tnorm{x}\le1}|f(x)|
 \qquad(f\in X^*).
\end{equation}
The star indicates a dual norm, not a derivative. In particular,
the two expressions in \eqref{eq:dual-norms-intro} need not have the
same value at the same functional. From
\eqref{eq:norm-equivalence-intro} one obtains
\[
 b^{-1}\norm{f}^{*}\le\tnorm{f}^{*}\le a^{-1}\norm{f}^{*}.
\]
Equivalent norms have the same continuous linear functionals, so
the underlying vector space \(X^*\) and the weak topology
\(\sigma(X,X^*)\) are unchanged. They also give equivalent operator
norms on spaces of continuous multilinear forms; these standard
consequences of norm equivalence are recorded in
\cite[Chapters~1--3]{FHHMZ}.

We write
\[
 B_{\tnorm{\cdot}}=\{x:\tnorm{x}\le1\},\qquad
 S_{\tnorm{\cdot}}=\{x:\tnorm{x}=1\}.
\]
Distances and open balls without an additional subscript refer to
\(\norm{\cdot}\). For derivatives of order at least two, an
unsubscripted operator norm is also induced by this original norm.
For first derivatives, which belong to \(X^*\), we use
\(\norm{\cdot}^{*}\). A functional norming \(x\ne0\) for the
equivalent norm satisfies
\[
 \tnorm{f}^{*}=1,\qquad f(x)=\tnorm{x}.
\]
Its original dual norm \(\norm{f}^{*}\) need not be one.

When \(X=\ell_p\), \(1<p<\infty\), the original norm is
\(\norm{\cdot}_p\). Put \(q=p/(p-1)\). The usual identification
\(X^*=\ell_q\) associates \(z=(z_k)\in\ell_q\) with
\(f_z(x)=\sum_k z_kx_k\), and
\[
 \norm{f_z}^{*}=\norm{z}_q,\qquad
 \tnorm{f_z}^{*}=\sup_{\tnorm{x}\le1}\left|\sum_k z_kx_k\right|.
\]
The duality identification and the first equality follow from
H\"older's inequality and its equality case; see
\cite[Chapter~1]{FHHMZ}. The second expression generally is a
different norm on \(\ell_q\).
Changing the background norm to an equivalent one does not change
whether a \emph{fixed function} is \Gat\ or \Fre\
differentiable. Replacing the norm function itself by
\(\tnorm{\cdot}\), however, may change both its differentiability
and the geometry of its unit ball; see \cite{Cudia,HJ}.

\subsection{First and higher derivatives}
\label{subsec:differentiability-intro}

Let \(U\subset X\) be open, let \(Y\) be a Banach space, and let
\(F:U\to Y\). The map \(F\) is \Gat\ differentiable at \(x\in U\)
if there is a bounded linear operator \(A:X\to Y\) such that
\begin{equation}\label{eq:G1-intro}
 \lim_{t\to0}\left\|
 \frac{F(x+th)-F(x)}{t}-Ah\right\|_Y=0
 \qquad\text{for every fixed }h\in X.
\end{equation}
We then write \(DF(x)=A\). Existence of all directional derivatives
alone is insufficient: their dependence on \(h\) must be linear and
continuous. The admissible size of \(t\) in
\eqref{eq:G1-intro} may depend on \(h\).

The map is \Fre\ differentiable at \(x\) if there is a bounded
linear \(A:X\to Y\) such that
\begin{equation}\label{eq:F1-intro}
 \lim_{\substack{u\to0\\u\ne0}}
 \frac{\norm{F(x+u)-F(x)-Au}_Y}{\norm{u}}=0.
\end{equation}
Equivalently, the convergence in \eqref{eq:G1-intro} is uniform
over \(\norm{h}=1\). Thus \Fre\ differentiability implies
\Gat\ differentiability, with the same derivative. The essential
extra requirement is uniformity over directions at the given base
point; it is not uniformity over all base points. These standard
notions and their renorming applications are treated in
\cite{HJ} and \cite[Chapters~7--8]{FHHMZ}.

For a norm, differentiability is always understood on
\(X\setminus\{0\}\). Indeed, for \(h\ne0\),
\[
 \frac{\tnorm{th}-\tnorm{0}}{t}
 =\frac{|t|}{t}\tnorm{h},
\]
which has no two-sided limit at zero. This obstruction is present
even for the Hilbert norm.

Higher \Fre\ derivatives are defined recursively by differentiating
the preceding derivative as a map into a space of bounded
multilinear forms equipped with its operator norm. The notation
\(C^m(U)\) additionally requires continuity of the derivatives up to
order \(m\). Mere existence of the \(m\)-th \Fre\ derivative at
every point should not silently be strengthened to continuity of
that derivative. The distinction matters even in finite dimensions; we use the
standard recursive calculus described in \cite{HJ}.

For a scalar-valued \(f\), the second \Fre\ derivative at \(x\)
is a bounded bilinear form \(A\) for which
\[
 \lim_{u\to0,\ u\ne0}
 \frac{\norm{Df(x+u)-Df(x)-A[u,\cdot]}^{*}}{\norm{u}}=0.
\]
At order three, a bounded trilinear form \(T\) must satisfy
\[
 \lim_{u\to0,\ u\ne0}
 \frac{\norm{D^2f(x+u)-D^2f(x)-T[u,\cdot,\cdot]}}{\norm{u}}=0,
\]
where the numerator now uses the operator norm of bilinear forms.
The preceding derivative must be defined in a neighborhood of
\(x\). For second or third \Gat\ differentiability in the
operator-norm convention, one tests the corresponding difference
quotients only along \(u=th\) for each fixed \(h\).

At higher orders we distinguish the following three \Gat\ notions.
First, this paper uses the \emph{recursive operator-norm convention}:
for each fixed \(h\), the difference quotient of \(D^{j-1}f\)
converges in the norm of the space of \((j-1)\)-linear forms to
\(D^jf(x)[h,\cdot,\ldots,\cdot]\). Definition~\ref{def:gateaux}
below states the formula precisely. This is the convention of
Maleev and Zlatanov \cite[Definition~2.1]{MaleevZlatanov}.
We denote the resulting class by \(G^m\).

Second, a weaker recursive convention only requires convergence
after evaluation at every fixed tuple of remaining arguments.
At order two, it asks for
\[
 \frac{Df(x+th)k-Df(x)k}{t}
 \longrightarrow D^2f(x)[h,k]
 \qquad(h,k\in X),
\]
with a bounded symmetric bilinear form on the right, rather than convergence
of the difference quotient in \(X^*\). We denote this convention
by \(G_{\mathrm w}^m\) when a distinction is needed.
It is used by McLaughlin and Vanderwerff
\cite[Section~1]{McLaughlinVanderwerff}, in their presentation of
Troyanski's higher-order theory.

Third, one may require only a \emph{directional polynomial
expansion} at each \(x\):
\begin{equation}\label{eq:directional-peano-intro}
 f(x+th)=f(x)+\sum_{j=1}^m\frac{t^j}{j!}A_j(x)[h^j]
                  +o_{x,h}(|t|^m),
\end{equation}
where each \(A_j(x)\) is a continuous symmetric \(j\)-linear form.
The remainder may depend on the fixed direction \(h\). This is the
definition used by Ivanov, Troyanski and Zlateva
\cite[Section~2]{ITZ}. Such an expansion does not by itself assert
recursive differentiability of the coefficient maps. We refer to
it as directional Peano differentiability.

The implications needed here are
\[
 \text{\(m\)-times \Fre}
 \ \Longrightarrow\ G^m
 \ \Longrightarrow\ G_{\mathrm w}^m
 \ \Longrightarrow\ \text{directional Peano of order \(m\)}.
\]
The first two arrows follow directly from the definitions; the
last is the ordinary one-variable Taylor--Peano formula applied
to each fixed affine line. Thus the arrows compare the conventions
of \cite{MaleevZlatanov,McLaughlinVanderwerff,ITZ}; they do not
identify them.
For example, at \(m=3\), a directional remainder \(o_{x,h}(|t|^3)\)
does not provide a common estimate for a sequence of different
directions; Theorem~\ref{thm:counter} gives an explicit example.
A theorem excluding the weakest class also excludes
the stronger classes; a construction in a weaker class needs
additional justification before it can be used in a stronger one.

\subsection{What smoothness says about convexity}
\label{subsec:geometry-intro}

Every norm is convex, independently of differentiability.
\emph{Strict convexity} means that its unit sphere contains no
nontrivial line segment. A norm is \emph{locally uniformly rotund}
(LUR) if, for \(x,x_n\in S_{\tnorm{\cdot}}\),
\[
 \tnorm{x+x_n}\longrightarrow2
 \quad\Longrightarrow\quad
 \tnorm{x_n-x}\longrightarrow0.
\]
Uniform convexity requires a corresponding implication uniformly
in both points. These are properties of the primal unit ball.
Smoothness has its most direct geometric interpretation through
functionals in the \emph{dual} unit ball
\cite{Day,Cudia}; see also \cite[Chapter~7]{FHHMZ}.

For \(x\ne0\), the norm \(\tnorm{\cdot}\) is \Gat\ differentiable
at \(x\) if and only if it has a unique norming functional
\cite{Cudia}, that is, a unique \(f_x\in X^*\) satisfying
\[
 \tnorm{f_x}^{*}=1,\qquad f_x(x)=\tnorm{x}.
\]
In that case \(D\tnorm{\cdot}(x)=f_x\). Thus smoothness means
uniqueness of the supporting functional at \(x\); it does not mean
that a supporting functional can touch the ball at only one point.
A flat portion of the sphere can have a unique supporting
functional at each of its points; Example~\ref{ex:smooth-flat}
provides a direct illustration.

For \Fre\ differentiability, the precise stronger statement is
the \v{S}mulyan criterion; see \cite{Cudia} and
\cite[Chapter~7]{FHHMZ}. If \(x\in S_{\tnorm{\cdot}}\),
then \(\tnorm{\cdot}\) is \Fre\ differentiable at \(x\) if and
only if there is a norming functional \(f_x\) such that
\begin{equation}\label{eq:smulyan-intro}
 \tnorm{f_n}^{*}\le1,\quad f_n(x)\longrightarrow1
 \quad\Longrightarrow\quad
 \tnorm{f_n-f_x}^{*}\longrightarrow0.
\end{equation}
Equivalently, the canonical image of \(x\) in \(X^{**}\) strongly
exposes the \emph{dual} unit ball at \(f_x\).
This does not say that \(x\) is strongly exposed in the
\emph{primal} unit ball.

The reason for the extra uniformity can be seen directly. If
\(\tnorm{\cdot}\) is \Fre\ differentiable at \(x\), then, for
small \(t>0\),
\[
 \tnorm{x\pm th}=1\pm t f_x(h)+t r_\pm(t,h),\qquad
 \sup_{\tnorm{h}\le1}|r_\pm(t,h)|\longrightarrow0.
\]
For \(\tnorm{f_n}^{*}\le1\), the inequalities
\(f_n(x\pm th)\le\tnorm{x\pm th}\) imply
\[
 |(f_n-f_x)(h)|
 \le\frac{1-f_n(x)}{t}
       +\max\{|r_+(t,h)|,|r_-(t,h)|\}.
\]
Taking the supremum in \(h\), then \(n\to\infty\), and finally
\(t\downarrow0\), gives \eqref{eq:smulyan-intro}.
Conversely, choose a functional \(g_u\) norming \(x+u\).
Then \(g_u(x)\ge1-2\tnorm{u}\), so
\eqref{eq:smulyan-intro} gives
\(\tnorm{g_u-f_x}^{*}\to0\). The supporting inequalities yield
\[
 0\le\tnorm{x+u}-1-f_x(u)
 \le (g_u-f_x)(u)
 \le\tnorm{g_u-f_x}^{*}\tnorm{u},
\]
which is the \Fre\ remainder estimate. This verifies the stated
version of the criterion directly.

Useful sufficient conditions run from the dual geometry to primal
smoothness \cite{Day,Cudia}: strict convexity of \(\tnorm{\cdot}^{*}\) implies
\Gat\ differentiability of \(\tnorm{\cdot}\), and LUR of
\(\tnorm{\cdot}^{*}\) implies \Fre\ differentiability of
\(\tnorm{\cdot}\). In a reflexive space, primal \Gat\
smoothness also implies strict convexity of the dual norm.
Indeed, if two distinct dual unit functionals had a midpoint of
norm one, that midpoint would attain its norm on the primal unit
ball; both functionals would then norm the same primal point,
contradicting uniqueness. Reflexivity is needed for this
norm-attainment argument. In general, one must not replace the
local dual statement \eqref{eq:smulyan-intro} by a claim that the
given primal norm is LUR.

De~Bernardi and Somaglia constructed, on every infinite-dimensional
separable Banach space, a strictly convex, \Gat\ smooth equivalent norm
with the Kadec property that fails to be LUR
\cite[Theorems~2.3 and~3.8]{DeBernardiSomaglia}.
De~Bernardi, Preti and Somaglia obtained a strictly convex, \Gat\
smooth norm for which even midpoint local uniform rotundity fails.
When the dual is separable, they obtained an equivalent norm that is
\Fre\ smooth and weakly uniformly rotund and still fails the midpoint
condition \cite[Theorems~3.3 and~3.4]{DeBernardiPretiSomaglia}.

Neither twice \Fre\ nor twice \Gat\ differentiability changes
this conclusion. If the second derivative of the norm exists,
convexity gives
\begin{equation}\label{eq:hessian-geometry-intro}
 D^2\tnorm{\cdot}(x)[h,h]\ge0.
\end{equation}
This is the usual second-derivative test for a convex line
restriction, within the calculus of \cite{HJ}; we give the
homogeneity calculation explicitly.
Positive homogeneity gives
\[
 D\tnorm{\cdot}(x)x=\tnorm{x},\qquad
 D^2\tnorm{\cdot}(x)[x,h]=0.
\]
The latter identity follows by differentiating the former
identity in direction \(h\). Thus the Hessian is necessarily
degenerate in the radial direction. Even on the tangent space
\(\ker D\tnorm{\cdot}(x)\), its nonnegativity does not imply
strict positivity or a uniform lower bound.
For example, if \(p>2\), then
\[
 D^2\norm{\cdot}_p(e_1)[e_2,e_2]=0.
\]
The vanishing follows by differentiating
\((1+|t|^p)^{1/p}\) twice at zero, whereas uniform convexity is
classical \cite{Hanner}.

\begin{example}[Smoothness does not force strict convexity]
\label{ex:smooth-flat}
A norm can be \(C^\infty\) away from zero and still have a flat
segment on its sphere. The following elementary construction
verifies this directly; the differentiability step uses the
implicit function theorem from Banach-space calculus \cite{HJ}.
Write \(\ell_2=\R^2\oplus_2 Z\).
Define the even, convex \(C^\infty\) function
\[
 \phi(s)=
 \begin{cases}
  0,&|s|\le1,\\[1mm]
  \displaystyle\int_1^{|s|}(|s|-u)
             \exp\!\left(-\frac{1}{(u-1)^2}\right)\dd u,
       &|s|>1.
 \end{cases}
\]
It is strictly increasing for \(s>1\) and tends to infinity.
The set
\[
 K=\{(s,t,z):\phi(s)+\phi(t)+\norm{z}_2^2\le1\}
\]
is closed, convex, symmetric and bounded, and contains a
neighborhood of zero. Its Minkowski functional
\(\tnorm{\cdot}_{\mathrm{flat}}\) is therefore an equivalent norm.
For a nonzero \((s,t,z)\), its value is the unique \(r>0\)
satisfying
\[
 \phi(s/r)+\phi(t/r)+\norm{z}_2^2/r^2=1.
\]
The derivative of the left side with respect to \(r\) is strictly
negative at a solution: at least one nonzero summand is active,
and each active summand is strictly decreasing there.
The implicit function theorem shows that this norm is \(C^\infty\)
off zero. If \(c>1\) is chosen with \(\phi(c)=1\), then
\((c,t,0)\), \(-1\le t\le1\), is a segment of its unit sphere.
The norm is consequently neither strictly convex nor LUR.
This example rules out any implication from first, second, third,
or even infinite-order smoothness to LUR of the \emph{same} norm.
\end{example}

\subsection{Existence results on general Banach spaces}
\label{subsec:general-survey}

Statements about the existence of a smooth equivalent norm are
isomorphic statements about \(X\); they do not describe every
equivalent norm on \(X\). At first order, every separable Banach
space admits an equivalent \Gat\ differentiable norm
\cite{Day}; an explicit construction, even with uniform
\Gat\ differentiability, is given in \cite[Section~1]{HMZ}.
For a separable Banach space, existence of an equivalent \Fre\
differentiable norm is equivalent to separability of \(X^*\)
\cite[Chapter~8]{FHHMZ}.
In arbitrary density, existence of a \Fre\ smooth equivalent
norm implies that \(X\) is Asplund, meaning that every separable
subspace has separable dual; no general converse is asserted here.
See Asplund \cite{Asplund} and
\cite[Chapter~11]{FHHMZ}. In particular, combining the preceding
separable renorming theorem with \((\ell_1)^*=\ell_\infty\) shows
that \Gat\ renormability does not imply the Asplund property.
Likewise, \(c_0\) has separable dual \(\ell_1\), so the same
\Fre\ renorming theorem gives a nonreflexive example; these
duality facts are recalled in \cite[Chapter~1]{FHHMZ}.

Asplund's averaging method gives every Banach space with separable
dual an equivalent norm that is both \(C^1\)-smooth and LUR
\cite{AsplundAveraged}. At higher orders, H\'ajek and Proch\'azka
\cite{HajekProchazka} showed that a weakly compactly generated Banach
space admitting a \(C^k\)-smooth norm, \(k\ge1\), has a \(C^1\)-smooth
LUR norm that is a uniform limit on bounded sets of \(C^k\)-smooth
norms. Their theorem gives first-order smoothness for the LUR norm
itself and order \(k\) for the approximating norms.

Higher-order existence results are often formulated for bump
functions. A bump is a nonzero continuous function with bounded
support. An equivalent norm of any of the smoothness classes used
here produces a bump of the same class: compose it with a smooth
scalar cutoff that is constant near zero and vanishes outside a
bounded interval. The composition is smooth at zero because it is
locally constant there, and its support is bounded by norm
equivalence. Consequently, nonexistence results for bumps imply
nonexistence results for equivalent norms. This standard passage from norms to bumps is used in
\cite{McLaughlinVanderwerff,DGJ}.

Vanderwerff \cite{Vanderwerff1993} proved that existence of a
continuous twice \Gat\ differentiable bump implies that the
space is Asplund. This obstruction applies already to the scalar
recursive convention, hence also to \(G^2\) and twice \Fre\
smoothness. A stronger conclusion holds with an additional
geometric hypothesis. McLaughlin, Poliquin, Vanderwerff and Zizler
\cite[Theorem~2.1 and Proposition~2.3]{MPVZ} showed that if \(X\)
has the Radon--Nikod\'ym property and admits such a bump, then it
admits an equivalent norm with modulus of smoothness of power
type two. For the latter norm this means
\[
 \rho_{\tnorm{\cdot}}(t):=
 \sup_{\tnorm{x}=\tnorm{y}=1}
 \left(\frac{\tnorm{x+ty}+\tnorm{x-ty}}2-1\right)\le Ct^2
 \quad(0\le t\le1)
\]
for some \(C\). Such a modulus implies superreflexivity and hence the existence
of an equivalent uniformly convex norm
\cite[Theorem~2 and its discussion]{HMZ}.
All reflexive spaces have the Radon--Nikod\'ym property
\cite[Chapter~11]{FHHMZ}.
The conclusion concerns a suitable equivalent norm, not
necessarily the original twice differentiable one.

For clarity, the geometric form of the Radon--Nikod\'ym property
used in this theory is dentability \cite[Section~1]{HMZ}:
every nonempty bounded
subset \(A\subset X\) has slices of arbitrarily small diameter.
A slice has the form
\[
 \{x\in A:f(x)>\sup_{y\in A}f(y)-\alpha\},
 \qquad f\in X^*,\quad\alpha>0.
\]
Thus this hypothesis supplies small slices of bounded sets; it
is an additional assumption on the space, not a definition of
smoothness.

An especially rigid result requires smoothness on both sides of
duality: if both \(X\) and \(X^*\) admit continuous twice \Gat\
differentiable bumps, then \(X\) is isomorphic to a Hilbert space
\cite{Vanderwerff1993}; see also \cite[Corollary~2.7 and the
following remark]{MPVZ}. Thus twice differentiability on \(X\)
alone and twice differentiability on both \(X\) and \(X^*\) must
not be confused.

Twice \Fre\ differentiability is a stronger requirement than
twice scalar-recursive \Gat\ differentiability. A classical
separation is provided by
\[
 X_0=\left(\bigoplus_{n=1}^{\infty}\ell_4^n\right)_{\ell_2}.
\]
This space admits an equivalent \(G_{\mathrm w}^2\)-smooth norm,
but no twice \Fre\ differentiable bump
\cite{MPVZ,DGJ}. At third order, \(\ell_3\) provides another
separation: it admits a third-order \Gat\ smooth equivalent
norm, including in the operator-norm convention discussed below,
but no three-times \Fre\ smooth equivalent norm
\cite{Troyanski,MaleevZlatanov,DGJ}.
There is no universal rule that second or third differentiability
of a norm off zero forces the given norm to be LUR, or forces the
space to be Hilbertian. For example, the canonical norm of
\(\ell_4\) is \(C^\infty\) off zero \cite{BonicFrampton},
and \(c_0\) also admits a \(C^\infty\) equivalent norm
\cite[Section~1]{HMZ}; \(c_0\) is not reflexive.
Quilis \cite{Quilis} constructed a \(C^\infty\)-smooth equivalent norm
on \(c_0\) whose closed unit ball is dentable and has no extreme points.
Example~\ref{ex:smooth-flat} addresses the geometry of the given
norm independently of these isomorphic existence results.

There is, however, a simple Hilbertian criterion involving the
origin. It is an elementary consequence of the second-order
Taylor formula \cite{HJ}, as the following calculation shows.
Suppose \(q(x)=\tnorm{x}^{\,2}\) is twice \Fre\
differentiable at zero. Its first derivative there is zero.
Its second-order expansion and \(q(th)=t^2q(h)\) give
\[
 t^2\tnorm{h}^{\,2}
 =\tfrac12t^2D^2q(0)[h,h]+o(t^2),
 \qquad
 \tnorm{h}^{\,2}=\tfrac12D^2q(0)[h,h].
\]
The form \(\tfrac12D^2q(0)\) is therefore an inner product
inducing \(\tnorm{\cdot}\). Thus this particular norm is a
Hilbert norm, and \(X\) is isomorphic to a Hilbert space.
The same reasoning works if the square has a second-order
directional expansion at zero with a single bounded bilinear
coefficient. This hypothesis is much stronger than twice
differentiability of the norm merely off zero.

\subsection{The first three orders on \texorpdfstring{\(\ell_p\)}{ell-p}}
\label{subsec:lp-survey}

We now consider \(\ell_p=\ell_p(\N)\), with \(1\le p<\infty\).
For \(1<p<\infty\), its canonical norm is uniformly convex
\cite{Hanner}. Uniform smoothness follows by applying uniform
convexity to the conjugate exponent and using smoothness--convexity
duality \cite[Chapter~9]{FHHMZ}. Consequently it is LUR and
\Fre\ differentiable away from zero. The standard differentiation
formulas underlying \cite{BonicFrampton} give
\begin{equation}\label{eq:lp-first-intro}
 D\norm{\cdot}_p(x)h
 =\norm{x}_p^{\,1-p}\sum_k |x_k|^{p-2}x_kh_k,\qquad x\ne0,
\end{equation}
where \(|x_k|^{p-2}x_k\) is interpreted as zero at \(x_k=0\).
H\"older's inequality shows that the coefficient sequence belongs
to \(\ell_q\) and has dual norm one.

For \(p>2\), the second derivative is
\begin{equation}\label{eq:lp-second-intro}
 \begin{split}
 D^2\norm{\cdot}_p(x)[h,k]
 &=(p-1)\norm{x}_p^{\,1-p}
       \sum_i |x_i|^{p-2}h_i k_i\\
 &\quad -(p-1)\norm{x}_p^{\,1-2p}
       \left(\sum_i |x_i|^{p-2}x_i h_i\right)
       \left(\sum_i |x_i|^{p-2}x_i k_i\right).
 \end{split}
\end{equation}
The first term is bounded by H\"older's inequality with exponents
\(p/(p-2),p,p\), and the second is a product of bounded linear
functionals. The formulas define continuous derivative maps on
\(\ell_p\setminus\{0\}\); thus the norm is \(C^2\) there.
More generally, for every positive integer \(m<p\), the canonical
norm is \(C^m\) off zero. If \(p\) is an even integer, then
\(\sum_k x_k^p\) is a continuous homogeneous polynomial and its
positive \(p\)-th root is real analytic off zero. These
classical smoothness conclusions are established in
\cite{BonicFrampton}; see also \cite[Introduction and
Section~II]{DGJ}.

For \(p\) not an even integer, the restriction
\[
 \norm{e_1+t e_2}_p=(1+|t|^p)^{1/p}
 =1+\frac1p|t|^p+O(|t|^{2p})
\]
shows the failure of the next integer order for the canonical
norm. If \(p\) is nonintegral, that order is \(\lceil p\rceil\);
if \(p\) is odd, it is \(p\). This one-dimensional calculation, included here explicitly,
concerns the canonical norm only; compare the sharp classical
orders in \cite{BonicFrampton}. To rule out \emph{every}
equivalent norm requires an isomorphic obstruction.

The classical \Fre\ obstruction is available: for a non-even
\(p\in(1,\infty)\), there is no \(m\)-times \Fre\ differentiable
equivalent norm when the integer \(m\ge p\)
\cite{DGJ,BonicFrampton}. In particular, no twice \Fre\ smooth
renorming exists for \(1<p<2\), and no three-times \Fre\ smooth
renorming exists for \(2<p\le3\).
The stronger bump obstruction in \cite[Corollary~1.2]{DGJ}
concerns \(T^p\)-smoothness: at each \(x\) there is a continuous
polynomial \(P_x\) of degree at most \(\lfloor p\rfloor\), with
\(P_x(0)=0\), such that
\[
 f(x+h)-f(x)-P_x(h)=o(\norm{h}^p)\qquad(h\to0).
\] An \(m\)-times \Fre\ smooth norm with \(m\ge p\)
would produce a bump with that property.

In the \Gat\ category, Fabian, Whitfield and Zizler
\cite{FWZ} exclude twice differentiable equivalent norms on
\(\ell_p\) for \(1\le p<2\); the corresponding bump obstruction
is discussed in \cite{McLaughlinVanderwerff,MPVZ}.
For odd integers, Troyanski \cite{Troyanski} obtained
equivalent \(p\)-times \Gat\ differentiable norms on
\(\sigma\)-finite \(L_p\)-spaces.
For the operator-norm convention specifically, the construction of
Maleev and Zlatanov \cite[Theorems~2 and~6]{MaleevZlatanov}
gives the positive result on \(\ell_3\): take the constant Nakano
exponent \(p_n=3\) in their Theorem~2.
Their conclusion is even uniform \Gat\ differentiability,
where the direction is fixed and the base point varies on the
sphere. This uniformity is different from uniformity over all
directions in \Fre\ differentiability.

Table~\ref{tab:orders} summarizes existence of an equivalent norm
at the first three orders. Every entry concerns all nonzero
points. The symbol \(F^m\) denotes existence of \(m\) recursive
\Fre\ derivatives, and \(G^m\) has the operator-norm meaning
fixed above. \emph{Unresolved here} means that neither the cited
results nor this paper settle unrestricted existence in that cell.
\begin{table}[htbp]
\centering
\small
\setlength{\tabcolsep}{4pt}
\begin{tabular}{@{}lccccccl@{}}
\hline
Exponent & \(G^1\) & \(F^1\) & \(G^2\) & \(F^2\) & \(G^3\) & \(F^3\) & Sources\\
\hline
\(p=1\)    & Yes & No  & No  & No  & No & No & \cite{Day,FHHMZ,FWZ}\\
\(1<p<2\)  & Yes & Yes & No  & No  & No & No & \cite{BonicFrampton,FWZ}\\
\(p=2\)    & Yes & Yes & Yes & Yes & Yes & Yes & \cite{BonicFrampton}\\
\(2<p<3\)  & Yes & Yes & Yes & Yes & Unresolved here & No & \cite{BonicFrampton,DGJ,ITZ}\\
\(p=3\)    & Yes & Yes & Yes & Yes & Yes & No & \cite{MaleevZlatanov,DGJ}\\
\(p>3\)    & Yes & Yes & Yes & Yes & Yes & Yes & \cite{BonicFrampton}\\
\hline
\end{tabular}
\caption{Existence of equivalent norms on the separable space
\(\ell_p\). The positive conclusions extend to every order when
\(p\) is an even integer. The \(G^3\) entry at \(p=3\) uses an equivalent
norm, not the canonical norm.}
\label{tab:orders}
\end{table}

For \(p=1\), the positive \(G^1\) entry follows from separability,
whereas \(F^1\) is impossible because the dual \(\ell_\infty\)
is not separable \cite[Chapters~1 and~8]{FHHMZ}.
For \(p>3\), the canonical norm already gives
all the positive entries through order three; the table makes
no assertion of infinite-order smoothness for a non-even \(p\)
\cite{BonicFrampton}.

The interval \(2<p<3\) is the focus of this paper. The canonical
norm is \(C^2\), but is not \(G^3\)-smooth on
\(\ell_p\setminus\{0\}\). This does not rule out another equivalent norm.
Ivanov, Troyanski and Zlateva \cite[Introduction]{ITZ} discuss
the unresolved exact order of \Gat\ smoothness in the
separable \(\ell_p(\N)\) setting. Their definition in
\cite[Section~2]{ITZ} uses the expansion \eqref{eq:directional-peano-intro}. Their Theorem~1.1 concerns
\(h_M(\Gamma)\) with \(\Gamma\) uncountable. For \(M(t)=|t|^p\)
and \(2<p<3\), it excludes a third-order directional Peano bump
on \(\ell_p(\Gamma)\), hence also a third-order smooth norm in
either recursive sense. The uncountability hypothesis cannot be
omitted when citing that conclusion for \(\ell_p(\N)\).
The older nonseparable odd-integer bump obstruction is given in
\cite{McLaughlinVanderwerff}. None of these statements supplies
an unrestricted nonexistence proof for the separable problem
addressed here.

\subsection{The Hilbert case \texorpdfstring{\(p=2\)}{p=2}}
\label{subsec:hilbert-intro}

The case \(\ell_2\) has no finite-order smoothness obstruction
\cite{BonicFrampton}. The following formulas are direct
derivatives of \(x\mapsto\sqrt{\langle x,x\rangle}\), using the
ordinary Banach-space chain rule \cite{HJ}.
Writing \(r=\norm{x}_2>0\), the canonical norm satisfies
\[
 \begin{split}
 D\norm{\cdot}_2(x)h
 &=\frac{\langle x,h\rangle}{r},\\
 D^2\norm{\cdot}_2(x)[h,k]
 &=\frac{\langle h,k\rangle}{r}
       -\frac{\langle x,h\rangle\langle x,k\rangle}{r^3},\\
 D^3\norm{\cdot}_2(x)[h,k,l]
 &=-\frac{\langle x,h\rangle\langle k,l\rangle
          +\langle x,k\rangle\langle h,l\rangle
          +\langle x,l\rangle\langle h,k\rangle}{r^3}
 +\frac{3\langle x,h\rangle\langle x,k\rangle
                      \langle x,l\rangle}{r^5}.
 \end{split}
\]
It is real analytic on \(\ell_2\setminus\{0\}\), since it is
the square root of the positive quadratic polynomial
\(\langle x,x\rangle\). Thus it is \(C^\infty\) there and has
\Gat\ derivatives of every order under all conventions above.
The squared norm is a quadratic polynomial on the whole space,
including zero:
\[
 D(\norm{\cdot}_2^2)(x)h=2\langle x,h\rangle,\qquad
 D^2(\norm{\cdot}_2^2)(x)[h,k]=2\langle h,k\rangle,\qquad
 D^j(\norm{\cdot}_2^2)=0\quad(j\ge3).
\]

The canonical norm is uniformly convex and uniformly smooth
\cite{Hanner}; the duality principle is recalled in
\cite[Chapter~9]{FHHMZ}. Its moduli, computed directly from
the parallelogram identity, are
\[
 \delta_{\ell_2}(\varepsilon)
 =1-\sqrt{1-\varepsilon^2/4}\quad(0\le\varepsilon\le2),
 \qquad
 \rho_{\ell_2}(t)=\sqrt{1+t^2}-1.
\]
Here \(\delta(\varepsilon)\) is the infimum of
\(1-\norm{(x+y)/2}_2\) over unit vectors at distance at least
\(\varepsilon\), and \(\rho\) is the modulus defined above.
The formulas follow from the parallelogram identity, with equality
realized in a two-dimensional subspace.
In particular, every point of the canonical sphere is an LUR
point and a strongly exposed point. For a unit vector \(x\), the
functional \(f_x(y)=\langle y,x\rangle\) strongly exposes \(x\),
because, for \(\norm{y}_2\le1\),
\[
 \norm{y-x}_2^2
 =\norm{y}_2^2+1-2\langle y,x\rangle
 \le2\bigl(1-f_x(y)\bigr).
\]
Thus \(f_x(y_n)\to1\) forces \(\norm{y_n-x}_2\to0\).
These statements describe the canonical Hilbert norm.
An arbitrary equivalent norm on \(\ell_2\) need not be smooth
or LUR, as the flat \(C^\infty\) example above already shows
for the latter property (Example~\ref{ex:smooth-flat}).

\subsection{Scope and organization of the present results}
\label{subsec:scope-intro}

Theorem~\ref{thm:baire} supplies the analytic step in the proof of
Theorem~\ref{thm:kadec}. A Baire argument in the product of base
points and directions gives a dense open set on which the third
derivative is locally bounded. The resulting cubic Taylor estimate
is uniform in the direction. At a strongly exposed point in this set,
Lemma~\ref{lem:cap} localizes a bump function. A related first-order
localization appears in \cite[Theorem~3.7]{BGV}. For a Kadec norm,
adding a compact Hilbertian term preserves the required smoothness
and makes every point of the new unit sphere strongly exposed
(Lemma~\ref{lem:kadec}).

The Baire-category and bump-function arguments have precedents in
\cite[Section~1]{HMZ} and \cite{McLaughlinVanderwerff,DGJ}.
In particular, \cite[Corollary~1.2]{DGJ} already rules out a nonzero
\(G^3\)-smooth bump function with globally bounded third derivative,
since a quadratic Taylor remainder of order \(O(\norm{h}_p^3)\) is
\(o(\norm{h}_p^p)\) for \(2<p<3\). The block argument in
Lemma~\ref{lem:bump} requires a uniform upper estimate only at points
where the bump function is at least \(1/2\). Theorem~\ref{thm:counter} shows precisely what
can fail at a single base point, even for a convex function.

Section~\ref{sec:prelim} fixes the definitions, and
Section~\ref{sec:main} states the three main results.
Sections~\ref{sec:aux} and~\ref{sec:proofs} contain the auxiliary
results and the proofs of the theorems, respectively.
Sections~\ref{sec:cor} and~\ref{sec:examples} give corollaries and examples.

\section{Preliminaries}\label{sec:prelim}

All Banach spaces and functions are real. For a Banach space \(X\),
let \(\Ls^j(X;\R)\) denote the continuous symmetric \(j\)-linear forms
on \(X\), with norm
\[
 \norm{A}=\sup\{\abs{A(h_1,\ldots,h_j)}:\norm{h_i}\le1,\ 1\le i\le j\}.
\]
We write \(A[h^j]=A(h,\ldots,h)\).
When \(j=1\), this is the dual norm \(\norm{A}^{*}\).
The symbols \(\norm{\cdot}^{*}\) and \(\tnorm{\cdot}^{*}\)
keep the meanings fixed in \eqref{eq:dual-norms-intro}.
Unless explicitly indicated otherwise, norms of vectors and forms on
\(\ell_p\) are taken with respect to the canonical norm
\[
 \norm{x}_p=\left(\sum_{k=1}^{\infty}|x_k|^p\right)^{1/p}.
\]
The standard unit vectors are denoted by \(e_k\).

\begin{definition}\label{def:gateaux}
Let \(U\subset X\) be open. A map \(f:U\to\R\) is \Gat\
differentiable at \(x\) if there exists \(Df(x)\in X^*\) such that
\[
 \lim_{t\to0}\frac{f(x+th)-f(x)}{t}=Df(x)h
 \qquad(h\in X).
\]
Inductively, suppose that \(D^{j-1}f\) is defined in a neighborhood of
\(x\). The \(j\)-th derivative is a form
\(D^jf(x)\in\Ls^j(X;\R)\) satisfying
\begin{equation}\label{eq:recursive}
 \lim_{t\to0}
 \left\|
 \frac{D^{j-1}f(x+th)-D^{j-1}f(x)}{t}
 -D^jf(x)[h,\cdot,\ldots,\cdot]
 \right\|_{\Ls^{j-1}}=0
 \qquad(h\in X).
\end{equation}
For brevity, we write \(Df(z)\) and \(D^2f(z)\) for the forms
\(Df(z)[\cdot]\) and \(D^2f(z)[\cdot,\cdot]\), respectively.
When a direction \(d\) is fixed, the dots in \(D^2f(z)[d,\cdot]\)
and \(D^3f(z)[d,\cdot,\cdot]\) indicate the arguments left free.

In particular, if \(Df\) exists in a neighborhood of \(z\), the
second \Gat\ derivative at \(z\) is the form \(D^2f(z)\in\Ls^2(X;\R)\)
satisfying
\[
 \lim_{t\to0}
 \left\|\frac{Df(z+td)-Df(z)}{t}-D^2f(z)[d,\cdot]\right\|^{*}=0
 \qquad(d\in X).
\]
Here \(D^2f(z)[d,\cdot]\) is the functional
\(v\mapsto D^2f(z)[d,v]\) in \(X^*\).
\end{definition}

\noindent\textbf{Smoothness classes.}
Let \(U\subset X\) be open and \(k\ge1\) an integer. For real-valued
functions on \(U\) we use the following notation.
\begin{description}[leftmargin=*,style=nextline]
\item \(G^k(U)\) be the recursive \Gat\ derivatives through order \(k\)
exist at each point of \(U\) in the sense of Definition~\ref{def:gateaux}.
\item \(F^k(U)\) be the recursive \Fre\ derivatives through order \(k\)
exist at each point of \(U\).
\item \(C^k(U)\) be the functions \(f\in F^k(U)\) for which \(D^kf\)
is continuous in operator norm on \(U\).
\item \(C^{2,1}_{\mathrm{loc}}(U)\) be the functions \(f\in C^2(U)\)
for which each \(x\in U\) has a neighborhood \(V\subset U\) and a
constant \(L\ge0\) satisfying
\(\norm{D^2f(y)-D^2f(z)}\le L\norm{y-z}\) for all \(y,z\in V\).
\end{description}
For a norm, these classes refer to \(U=X\setminus\{0\}\) unless
a smaller domain is specified.

Whenever \(D^3f(x)\) exists in the sense of Definition~\ref{def:gateaux},
it is a continuous symmetric trilinear form, so
\(\norm{D^3f(x)}<\infty\). Its diagonal \(h\mapsto D^3f(x)[h^3]\)
is the associated cubic polynomial. For \(f\in G^3(U)\), we say
that \(D^3f\) is \emph{locally bounded at \(x\)} if there are a
neighborhood \(V\subset U\) of \(x\) and a constant \(L<\infty\)
such that \(\norm{D^3f(y)}\le L\) for every \(y\in V\).
Thus the same bound must hold as the point \(y\) varies in \(V\).

The operator norm in \eqref{eq:recursive} is essential to this
convention. The arguments also use that the restrictions to affine
lines possess ordinary derivatives up to order three. In particular,
\[
 f(x+th)=f(x)+tDf(x)h+\frac{t^2}{2}D^2f(x)[h^2]
       +\frac{t^3}{6}D^3f(x)[h^3]+o_{x,h}(|t|^3).
\]
No uniformity in \(h\) is part of this formula.

A \emph{globally Lipschitz second
derivative} means that one Lipschitz constant works on the whole
space. A bump function is a nonzero continuous real-valued function
with bounded support, where
\(\supp b=\overline{\{x:b(x)\ne0\}}\).

\begin{definition}\label{def:geometry}
Let \(\norm{\cdot}\) be a norm on \(X\), and put
\(B_{\norm{\cdot}}=\{x:\norm{x}\le1\}\), \(S_{\norm{\cdot}}=\{x:\norm{x}=1\}\).
A point \(x_0\in S_{\norm{\cdot}}\) is strongly exposed if there exists
\(f\in X^*\), with \(\norm{f}^{*}=f(x_0)=1\), such that
\[
 y_n\in B_{\norm{\cdot}},\quad f(y_n)\to1
 \quad\Longrightarrow\quad \norm{y_n-x_0}\to0.
\]
The norm \(\norm{\cdot}\) has the Kadec property if its weak and norm
topologies agree on \(S_{\norm{\cdot}}\). On the separable reflexive spaces used
below, this is equivalently expressed by
\[
 x_n,x\in S_{\norm{\cdot}},\quad x_n\rightharpoonup x
 \quad\Longrightarrow\quad \norm{x_n-x}\to0.
\]
The norm is locally uniformly rotund (LUR) if, for every
\(x,x_n\in S_{\norm{\cdot}}\), the condition \(\norm{x+x_n}\to2\) implies
\(\norm{x_n-x}\to0\).
\end{definition}

These standard notions are discussed in \cite{FHHMZ}. In particular,
every LUR norm is Kadec. A norming functional at an LUR point
strongly exposes that point.

A normalized successive block sequence is a sequence
\((u_n)\) of finitely supported vectors with
\(\norm{u_n}_p=1\) and
\(\max\supp u_n<\min\supp u_{n+1}\).
The compact quadratic form used throughout is
\begin{equation}\label{eq:QM}
 Q(x)=\sum_{k=1}^{\infty}2^{-k}x_k^2,\qquad
 M(x)=\norm{x}_p^p.
\end{equation}
For \(p>2\), \(Q\) is a continuous quadratic polynomial,
\(Q(x)>0\) for \(x\ne0\), and \(Q(x)\le\norm{x}_p^2\).
The function \(M\) is \(C^2\), with
\[
 DM(x)h=p\sum_k |x_k|^{p-2}x_kh_k,\qquad
 D^2M(x)[h,k]=p(p-1)\sum_i |x_i|^{p-2}h_i k_i.
\]
The boundedness of these expressions follows from H\"older's
inequality; their continuity is the usual \(L_p\) differentiation
property.

\section{Main results}\label{sec:main}

Let \(U\) be an open subset of a real Banach space \(X\) and let
\(f\in G^3(U)\). We write
\[
 G_f(U)=\left\{x\in U:\ \exists\,r>0,\quad B(x,r)\subset U,\quad
          \sup_{y\in B(x,r)}\norm{D^3f(y)}<\infty\right\}.
\]
Thus \(x\in G_f(U)\) means that \(D^3f\) is locally bounded at \(x\)
in operator norm: one finite bound applies to \(D^3f(y)\) for every
\(y\) in some neighborhood of \(x\). When \(f\) is a norm
\(\tnorm{\cdot}\) on \(X\), we take \(U=X\setminus\{0\}\) and write
\(G_{\tnorm{\cdot}}\) for \(G_{\tnorm{\cdot}}(X\setminus\{0\})\).

We first state the regularity result that supplies common Taylor
estimates. It applies to an arbitrary real Banach space.

\begin{theorem}[Generic local boundedness]\label{thm:baire}
Let \(U\subset X\) be open, and let \(f:U\to\R\) be continuous and
belong to \(G^3(U)\). Then the following assertions hold.
\begin{enumerate}[label=\textup{\arabic*.},ref=\arabic*,leftmargin=*]
\item\label{it:baire-open} \(G_f(U)\) is open and dense in \(U\).
\item The restriction of \(f\) to \(G_f(U)\) belongs to
\(C^{2,1}_{\mathrm{loc}}(G_f(U))\).
\item\label{it:baire-taylor} If \(B\subset U\) is an open ball on which
\(\norm{D^3f}\le L\), then
\begin{equation}\label{eq:uniform-taylor}
 \abs{f(x+h)-f(x)-Df(x)h-\tfrac12D^2f(x)[h^2]}
 \le \frac L6\norm{h}^3
\end{equation}
whenever the segment \([x,x+h]\) is contained in \(B\).
\end{enumerate}
\end{theorem}

The next theorem is the main geometric application. Its conclusion
is conditional on the existence of an unrestricted \(G^3\)-smooth
norm; it excludes additional geometric properties of that norm.

\begin{theorem}[Kadec obstruction]\label{thm:kadec}
Let \(2<p<3\). Suppose that \(\tnorm{\cdot}\) is an equivalent \(G^3\)-smooth
norm on \(\ell_p\). Then the following assertions hold.
\begin{enumerate}[label=\textup{\arabic*.},ref=\arabic*,leftmargin=*]
\item\label{it:kadec-regular} \(G_{\tnorm{\cdot}}\cap S_{\tnorm{\cdot}}\) is relatively open and dense in \(S_{\tnorm{\cdot}}\).
No point of \(G_{\tnorm{\cdot}}\cap S_{\tnorm{\cdot}}\) is a strongly exposed point of \(B_{\tnorm{\cdot}}\).
\item\label{it:kadec-exposed} The strongly exposed points of \(B_{\tnorm{\cdot}}\) are contained in a
relatively closed nowhere dense subset of \(S_{\tnorm{\cdot}}\).
\item The norm \(\tnorm{\cdot}\) does not have the Kadec property. Consequently,
\(\ell_p\) admits no equivalent norm that is both \(G^3\)-smooth and Kadec.
\end{enumerate}
\end{theorem}

The final theorem shows why a fixed-base-point argument cannot
replace the global analysis in Theorem~\ref{thm:kadec}.

\begin{theorem}[Failure of block-uniform cubic remainders]
\label{thm:counter}
Let \(2<p<3\), choose an integer \(m\ge2\) with \(m(p-2)>1\), and set
\begin{equation}\label{eq:F}
 \Psi(x)=\bigl(Q(x)^m+M(x)^m\bigr)^{1/m}.
\end{equation}
Then the following assertions hold.
\begin{enumerate}[label=\textup{\arabic*.},ref=\arabic*,leftmargin=*]
\item \(\Psi\) is continuous, even and convex, is \(C^1\) on \(\ell_p\),
and is \(C^2\) on \(\ell_p\setminus\{0\}\).
\item \(\Psi\) is twice \Fre\
differentiable everywhere and three times \Gat\ differentiable
at zero, with
\begin{equation}\label{eq:jets-F}
 D\Psi(0)=0,\qquad D^2\Psi(0)=D^2Q,\qquad D^3\Psi(0)=0.
\end{equation}
\item For every \(h\in\ell_p\),
\begin{equation}\label{eq:coercive-F}
 \Psi(h)-\Psi(0)-D\Psi(0)h\ge\norm{h}_p^p.
\end{equation}
\item For every normalized successive block sequence
\((u_n)\) and every \(t\ne0\),
\begin{equation}\label{eq:sup-remainder}
 \sup_n\bigl|\Psi(tu_n)-t^2Q(u_n)\bigr|=|t|^p.
\end{equation}
In particular, no such sequence, or block subsequence, admits a
common \(O(|t|^3)\) bound for these remainders near zero.
\end{enumerate}
\end{theorem}

In the notation of Section~\ref{sec:prelim},
\(\Psi\in F^2(\ell_p)\setminus C^2(\ell_p)\).
Example~\ref{ex:hessian-discontinuous} shows that \(D^2\Psi\)
is discontinuous in operator norm at zero.

\section{Auxiliary results}\label{sec:aux}

In an expression such as \(D^3f(z)[h,u,v]\), the \emph{base point}
\(z\) is the point at which the derivative is evaluated, and \(h,u,v\)
are the directions to which the resulting trilinear form is applied.
A bound that holds uniformly for all base points \(z\in B\) uses the
same constant for every \(z\in B\). Given an open ball \(B=B(a,R)\),
a \emph{smaller concentric ball} is an open ball \(B(a,r)\) with
the same center \(a\) and radius \(0<r<R\).

The first elementary observation supplies the directions in the
moving-coordinate argument. It is a standard diagonal-vanishing
phenomenon; the proof is included to avoid invoking an operator
compactness theorem.

\begin{lemma}[Vanishing on disjoint blocks]\label{lem:diagonal}
Let \(p>2\), let \(B\) be a continuous bilinear form on \(\ell_p\),
and let \((u_n)\) be a normalized disjointly supported sequence.
Then 
\begin{enumerate}
\item \(B(u_n,u_n)\to0\)
\item \(u_n\rightharpoonup0\).
\end{enumerate}
\end{lemma}
\begin{proof}

1.) We start with the proof of \(B(u_n,u_n)\to0\) and we will divide it into steps.

\smallskip
\noindent\emph{Step 1:  We will prove the inequality \(|B(v,v)|\le\norm{B}\norm{v}_p^2\).}

Write \(u_n=(u_n(k))_{k\ge1}\) and set
\( I_n=\supp u_n=\{k\in\N:u_n(k)\ne0\}\).
The assumptions mean that
\(I_n\cap I_m=\varnothing\) for \(n\ne m\) and \(\sum_{k\in I_n}|u_n(k)|^p=\norm{u_n}_p^p=1\).
Thus, at each coordinate \(k\), at most one vector of the sequence
has a nonzero entry. In particular, if \(n_1,\ldots,n_r\) are
distinct indices and \(a_1,\ldots,a_r\in\R\), then
\begin{equation}\label{eq:diag-disjoint-norm}
 \left\|\sum_{j=1}^{r}a_ju_{n_j}\right\|_p^p
 =\sum_{k=1}^{\infty}
       \left|\sum_{j=1}^{r}a_ju_{n_j}(k)\right|^p
 =\sum_{j=1}^{r}|a_j|^p
       \sum_{k\in I_{n_j}}|u_{n_j}(k)|^p
 =\sum_{j=1}^{r}|a_j|^p.
\end{equation}

There is no inequality in this calculation: disjointness prevents
different vectors from contributing to the same coordinate, and
normalization makes each inner sum equal to one. Consequently, for
every choice \(\sigma_j\in\{-1,1\}\),
\begin{equation}\label{eq:diag-signed-norm}
 \left\|\sum_{j=1}^{r}\sigma_ju_{n_j}\right\|_p=r^{1/p}.
\end{equation}

The supports need not be finite, and their order in \(\N\) is
irrelevant. Pairwise disjointness is all that is used here. A
normalized successive block sequence, as defined in the
preliminaries, is therefore a particular case of the present
assumptions.

Continuity and bilinearity of \(B\) imply that its operator norm
\[
 \norm{B}=\sup\{|B(v,w)|:\norm{v}_p\le1,\ \norm{w}_p\le1\}
\]
is finite. Rescaling each nonzero argument gives
\begin{equation}\label{eq:diag-bilinear-bound}
 |B(v,w)|\le\norm{B}\,\norm{v}_p\norm{w}_p
 \qquad(v,w\in\ell_p).
\end{equation}

If either argument is zero, the same inequality follows from
bilinearity. In particular,
\(|B(v,v)|\le\norm{B}\norm{v}_p^2\).

\smallskip
\noindent\emph{Step 2: We will prove that, if \(B(u_n,u_n)\not\to0\),
there are \(\varepsilon>0\) and a subsequence \((u_{n_j})\) such that
\(\bigl|\sum_{j=1}^{r}B(u_{n_j},u_{n_j})\bigr|\ge r\varepsilon\)
for every \(r\in\N\).}

Suppose, towards a contradiction, that \(B(u_n,u_n)\) does not
converge to zero. The negation of convergence gives an
\(\varepsilon>0\) for which
\( \{n\in\N:|B(u_n,u_n)|\ge\varepsilon\}\)
is infinite. Since the scalar field is real, at least one of the
sets
\(\{n:B(u_n,u_n)\ge\varepsilon\}\) or \(\{n:B(u_n,u_n)\le-\varepsilon\}\)
is infinite. We can therefore select increasing indices
\(n_1<n_2<\cdots\) such that all the corresponding diagonal
values have the same sign and modulus at least \(\varepsilon\).
It follows that, for every positive integer \(r\),
\begin{equation}\label{eq:diag-lower-bound}
 \left|\sum_{j=1}^{r}B(u_{n_j},u_{n_j})\right|
 =\sum_{j=1}^{r}|B(u_{n_j},u_{n_j})|
 \ge r\varepsilon.
\end{equation}

\smallskip
\noindent\emph{Step 3: We will prove the identity
\(2^{-r}\sum_{\sigma\in\{-1,1\}^r}B(v_\sigma,v_\sigma)
=\sum_{j=1}^{r}B(u_{n_j},u_{n_j})\), where
\(v_\sigma=\sum_{j=1}^{r}\sigma_ju_{n_j}\).}

Fix \(r\in\N\). Let
\( \Sigma_r=\{-1,1\}^r\) and \(v_\sigma=\sum_{j=1}^{r}\sigma_ju_{n_j}\) for \(\sigma=(\sigma_1,\ldots,\sigma_r)\in\Sigma_r\).
There are exactly \(2^r\) sign choices. We take the ordinary
arithmetic mean of the \(2^r\) real numbers \(B(v_\sigma,v_\sigma)\).
Expanding both arguments by bilinearity and interchanging finite
sums gives
\[
 \begin{split}
 2^{-r}\sum_{\sigma\in\Sigma_r}B(v_\sigma,v_\sigma)
 &=2^{-r}\sum_{\sigma\in\Sigma_r}
       \sum_{i=1}^{r}\sum_{j=1}^{r}
       \sigma_i\sigma_jB(u_{n_i},u_{n_j})\\
 &=\sum_{i=1}^{r}\sum_{j=1}^{r}
       \left(2^{-r}\sum_{\sigma\in\Sigma_r}\sigma_i\sigma_j\right)
       B(u_{n_i},u_{n_j}).
 \end{split}
\]

For \(i=j\), we have \(\sigma_i\sigma_j=\sigma_i^2=1\) for
every sign choice, so the coefficient in parentheses equals one.
For \(i\ne j\), pair each \(\sigma\) with the sign vector
obtained by reversing only its \(i\)-th sign. This pairing is a
bijection without fixed points. The value of \(\sigma_i\sigma_j\)
changes to its negative, because \(\sigma_j\) is unchanged.
Each pair therefore contributes zero. In symbols,
\[
 2^{-r}\sum_{\sigma\in\Sigma_r}\sigma_i\sigma_j
 =\begin{cases}
   1,&i=j,\\
   0,&i\ne j.
  \end{cases}
\]

Consequently,
\begin{equation}\label{eq:diag-averaging}
 2^{-r}\sum_{\sigma\in\Sigma_r}B(v_\sigma,v_\sigma)
 =\sum_{j=1}^{r}B(u_{n_j},u_{n_j}).
\end{equation}

\smallskip
\noindent\emph{Step 4: We will prove that \(B(u_n,u_n)\to0\)
by deriving the impossible bound
\(\varepsilon\le\norm{B}r^{2/p-1}\) for every \(r\in\N\).}

By \eqref{eq:diag-signed-norm} and
\eqref{eq:diag-bilinear-bound}, each summand satisfies
\( |B(v_\sigma,v_\sigma)|
 \le\norm{B}\norm{v_\sigma}_p^2
 =\norm{B}r^{2/p}\).
Combining the lower bound \eqref{eq:diag-lower-bound} with the
averaging identity \eqref{eq:diag-averaging} and the triangle
inequality now yields
\[
 \begin{split}
 r\varepsilon
 &\le\left|\sum_{j=1}^{r}B(u_{n_j},u_{n_j})\right|
 =\left|2^{-r}\sum_{\sigma\in\Sigma_r}B(v_\sigma,v_\sigma)\right|
 \le2^{-r}\sum_{\sigma\in\Sigma_r}|B(v_\sigma,v_\sigma)|\\
 &\le2^{-r}\,2^r\norm{B}r^{2/p}
 =\norm{B}r^{2/p}.
 \end{split}
\]

Dividing by \(r\) gives
\( \varepsilon\le\norm{B}r^{2/p-1}
 =\frac{\norm{B}}{r^{(p-2)/p}}\) for \(r\in\N\).
The exponent
\((p-2)/p\) is positive, so the right-hand side tends to zero as
\(r\to\infty\). The left-hand side is the fixed positive number
\(\varepsilon\), which is impossible. Hence
\( B(u_n,u_n)\longrightarrow0\).

2.) We will prove that \(f(u_n)\to0\)
for every \(f\in(\ell_p)^*\), and hence \(u_n\rightharpoonup0\).

Fix \(f\in\ell_p^*\), and put \(p'=p/(p-1)\), so that
\(1/p+1/p'=1\). The standard dual identification
\((\ell_p)^*=\ell_{p'}\) means that there is a sequence
\(a=(a_k)_{k\ge1}\in\ell_{p'}\) such that
\(f(x)=\sum_{k=1}^{\infty}a_kx_k\) for \(x\in\ell_p\) and \(\norm{f}^{*}=\norm{a}_{p'}\).

The series is absolutely convergent by H\"older's inequality.

For each support \(I_n\), define
\( c_n=\left(\sum_{k\in I_n}|a_k|^{p'}\right)^{1/p'}\).
Disjointness of the sets \(I_n\) gives, for every \(r\in\N\),
\[
 \sum_{n=1}^{r}c_n^{p'}
 =\sum_{k\in I_1\cup\cdots\cup I_r}|a_k|^{p'}
 \le\sum_{k=1}^{\infty}|a_k|^{p'}
 =\norm{a}_{p'}^{p'}<\infty.
\]

The nonnegative series \(\sum_n c_n^{p'}\) therefore converges,
so its terms tend to zero, and hence \(c_n\to0\). Equivalently,
if infinitely many \(c_n\) were at least some fixed
\(\eta>0\), the partial sums would contain arbitrarily many
terms of size at least \(\eta^{p'}\), contradicting their common
finite bound.

Since \(u_n(k)=0\) for \(k\notin I_n\), H\"older's inequality
on \(I_n\) and normalization give
\[
 |f(u_n)|
 =\left|\sum_{k\in I_n}a_ku_n(k)\right|
 \le\left(\sum_{k\in I_n}|a_k|^{p'}\right)^{1/p'}
       \left(\sum_{k\in I_n}|u_n(k)|^p\right)^{1/p}
 =c_n\norm{u_n}_p=c_n\longrightarrow0.
\]

This proves the second assertion. In particular, \((u_n)\) is
weakly null.
\end{proof}

\begin{lemma}[Bounded third derivatives]\label{lem:mvt}
Let \(f\in G^3(B)\), where \(B\) is a convex open ball, and suppose
\(\sup_B\norm{D^3f}\le L\). Then 
\begin{enumerate}[label=\textup{(\roman*)},leftmargin=*]
\item \(D^2f\) is Lipschitz on \(B\) with constant \(L\).
\item \(f\in C^{2,1}_{\mathrm{loc}}(B)\), and \eqref{eq:uniform-taylor} holds.
\item Both \(Df\) and \(D^2f\) are bounded on \(B\).
\end{enumerate}
\end{lemma}

\begin{proof}
The proof applies to an arbitrary real Banach space \(X\). Write
\(B=B(a,R)=\{x\in X:\norm{x-a}<R\}\), where \(a\in X\) and
\(R>0\). The norms of \(Df(x)\), \(D^2f(x)\), and \(D^3f(x)\)
are the corresponding operator norms. In particular, the hypothesis
means that
\begin{equation}\label{eq:mvt-third-bound}
 |D^3f(z)[h,u,v]|
 \le L\norm{h}\norm{u}\norm{v}
 \qquad(z\in B,\ h,u,v\in X).
\end{equation}
The same constant \(L\) works for every \(z\in B\) and all \(h,u,v\in X\).

We will prove norm continuity of \(Df\) and \(D^2f\) as functions
of the base point. The proof uses only the stated bound on \(D^3f\).
The order of the argument is therefore important: first we obtain
Lipschitz continuity of \(D^2f\) from the scalar mean value theorem;
only then do we use integration to identify the \Fre\ derivatives
and estimate the Taylor remainder.

\smallskip
\noindent\emph{Part (i): We will prove that
\(\norm{D^2f(y)-D^2f(x)}\le L\norm{y-x}\) for all
\(x,y\in B\).}

Fix \(x,y\in B\), and put \(d=y-x\). If \(d=0\), the desired
Lipschitz inequality is immediate, so suppose that \(d\ne0\).
Convexity ensures that \(x+sd\in B\) for \(0\le s\le1\).
Moreover, the set of real \(s\) for which \(x+sd\in B\) is open
and contains \([0,1]\). Thus the line restrictions below are
defined on an open interval containing this closed interval,
including small neighborhoods of its endpoints.

For fixed vectors \(u,v\in X\), consider the real-valued function
\( \varphi_{u,v}(s)=D^2f(x+sd)[u,v]\).
By the recursive definition of the third \Gat\ derivative, at
\(z=x+sd\) we have
\[
 \frac{D^2f(z+td)-D^2f(z)}{t}
 \longrightarrow D^3f(z)[d,\cdot,\cdot]
 \quad\text{in the operator norm as }t\to0.
\]

Evaluation at \((u,v)\) is continuous for that norm, since
\(|A[u,v]|\le\norm{A}\norm{u}\norm{v}\) for every bounded
bilinear form \(A\). For \(t\ne0\) sufficiently small, apply this
estimate to the bilinear form
\[
 A_t=\frac{D^2f(z+td)-D^2f(z)}{t}-D^3f(z)[d,\cdot,\cdot].
\]
The preceding limit gives \(\norm{A_t}\to0\), so
\(\abs{A_t[u,v]}\le\norm{A_t}\norm{u}\norm{v}\to0\).
Since
\[
 A_t[u,v]=\frac{\varphi_{u,v}(s+t)-\varphi_{u,v}(s)}{t}
           -D^3f(z)[d,u,v],
\]
we obtain
\[
 \varphi_{u,v}'(s)=D^3f(x+sd)[d,u,v].
\]

In particular, \(\varphi_{u,v}\) is an ordinary differentiable
real function. It is consequently continuous on \([0,1]\) and
differentiable on \((0,1)\), which are precisely the hypotheses
needed for the scalar mean value theorem. Continuity of
\(\varphi_{u,v}'\) is not required.

The mean value theorem and \eqref{eq:mvt-third-bound} give
\[
 |(D^2f(y)-D^2f(x))[u,v]|
 =|\varphi_{u,v}(1)-\varphi_{u,v}(0)|
 \le L\norm{y-x}\norm{u}\norm{v}.
\]

Although the intermediate point supplied by the mean value
theorem can depend on \(u\) and \(v\), the bound on the right
does not depend on that point. Taking the supremum over
\(\norm{u}\le1\) and \(\norm{v}\le1\) proves
\begin{equation}\label{eq:mvt-hessian-lipschitz}
 \norm{D^2f(y)-D^2f(x)}\le L\norm{y-x}
 \qquad(x,y\in B).
\end{equation}

\smallskip
\noindent\emph{Part (ii), Step 1: We will prove that
\(\norm{Df(x+h)-Df(x)}^{*}
\le\norm{D^2f(x)}\norm{h}+\tfrac L2\norm{h}^2\)
whenever \([x,x+h]\subset B\), and deduce that \(Df\) is
norm-continuous on \(B\).}

Fix \(x\in B\) and \(h\in X\) with \([x,x+h]\subset B\).
For a fixed \(z\in X\), define
\[
 \psi_z(s)=Df(x+sh)z,\qquad 0\le s\le1.
\]

The definition of \(D^2f\), followed by evaluation at \(z\), gives
\(\psi_z'(s)=D^2f(x+sh)[h,z]\).
This derivative is continuous in \(s\), because
\eqref{eq:mvt-hessian-lipschitz}, proved in part~(i), establishes
norm-continuity of \(D^2f\). The fundamental
theorem of calculus is therefore applicable and yields
\begin{equation}\label{eq:mvt-gradient-identity}
 (Df(x+h)-Df(x))z
 =\int_0^1D^2f(x+sh)[h,z]\dd s.
\end{equation}

All integrals in this argument can be understood as ordinary
real integrals after evaluation at a fixed vector. 

By \eqref{eq:mvt-hessian-lipschitz}, we obtain
\(\norm{D^2f(x+sh)}
 \le\norm{D^2f(x)}+Ls\norm{h}\).
Using this estimate in \eqref{eq:mvt-gradient-identity}, and then
taking the supremum over \(\norm{z}\le1\), gives
\begin{equation}\label{eq:mvt-gradient-increment}
 \norm{Df(x+h)-Df(x)}^{*}
 \le\norm{h}\int_0^1\norm{D^2f(x+sh)}\dd s
 \le\norm{D^2f(x)}\norm{h}+\frac L2\norm{h}^2.
\end{equation}

For each fixed \(x\), the right-hand side tends to zero as
\(\norm{h}\to0\). Hence \(Df\) is norm-continuous on \(B\).

For later use, we record a local Lipschitz bound as well. If
\(0<r<R\), then comparison with the center \(a\) in
\eqref{eq:mvt-hessian-lipschitz} gives
\(\norm{D^2f(y)}\le K_r:=\norm{D^2f(a)}+Lr\) for \(y\in B(a,r)\).

For \(x,y\in B(a,r)\), the segment joining them remains in the
same smaller ball. Formula~\eqref{eq:mvt-gradient-identity},
with \(h=y-x\), therefore implies
\begin{equation}\label{eq:mvt-gradient-local-lipschitz}
 \norm{Df(y)-Df(x)}^{*}\le K_r\norm{y-x}
 \qquad(x,y\in B(a,r)).
\end{equation}

\smallskip
\noindent\emph{Step 2: We will prove that \(f\) is continuously
\Fre\ differentiable on \(B\), with its \Fre\ derivative
equal to the given \Gat\ derivative \(Df\).}

We now verify the \Fre\ limit explicitly. Fix \(x\in B\), and
choose \(r\) with \(\norm{x-a}<r<R\). If \(\norm{h}\) is
sufficiently small, then \([x,x+h]\subset B(a,r)\).
The function \(s\mapsto f(x+sh)\) has ordinary derivative
\(Df(x+sh)h\), which is continuous by Step~1 of part~(ii). Hence
\[
 f(x+h)-f(x)=\int_0^1Df(x+sh)h\dd s.
\]

Subtracting \(Df(x)h\) and using
\eqref{eq:mvt-gradient-local-lipschitz} gives
\[
 \begin{split}
 |f(x+h)-f(x)-Df(x)h|
 &\le\int_0^1\norm{Df(x+sh)-Df(x)}^{*}\norm{h}\dd s\\
 &\le K_r\norm{h}^2\int_0^1s\dd s
 =\frac{K_r}{2}\norm{h}^2.
 \end{split}
\]

For \(h\ne0\), division by \(\norm{h}\) makes the right-hand
side equal to \(K_r\norm{h}/2\), which tends to zero.
The bound holds for all sufficiently small vectors \(h\), with
no restriction to a fixed direction. This is exactly \Fre\
differentiability of \(f\) at \(x\), and its \Fre\ derivative
is the previously given \Gat\ derivative \(Df(x)\).
Together with the norm-continuity of \(Df\) established in
Step~1 of part~(ii), this proves that \(f\) is continuously
\Fre\ differentiable on \(B\).

\smallskip
\noindent\emph{Step 3: We will prove that \(Df\) is \Fre\
differentiable at each \(x\in B\), with derivative
\(h\mapsto D^2f(x)[h,\cdot]\), and conclude that
\(f\in C^{2,1}_{\mathrm{loc}}(B)\).}

Fix \(x,h\) as in Step~2. For every \(z\in X\), subtract
\(D^2f(x)[h,z]\) from both sides of
\eqref{eq:mvt-gradient-identity}. Since this quantity is independent
of \(s\), the result is
\[
 (Df(x+h)-Df(x))z-D^2f(x)[h,z]
 \qquad=\int_0^1
       (D^2f(x+sh)-D^2f(x))[h,z]\dd s.
\]

The Lipschitz estimate \eqref{eq:mvt-hessian-lipschitz} implies
\[
 |(Df(x+h)-Df(x))z-D^2f(x)[h,z]|
 \le\int_0^1Ls\norm{h}^2\norm{z}\dd s
 =\frac L2\norm{h}^2\norm{z}.
\]

Taking the supremum over \(\norm{z}\le1\) yields the
operator-norm remainder estimate
\begin{equation}\label{eq:mvt-first-derivative-remainder}
 \norm{Df(x+h)-Df(x)-D^2f(x)[h,\cdot]}^{*}
 \le\frac L2\norm{h}^2.
\end{equation}

The map
\(h\longmapsto D^2f(x)[h,\cdot]\)
is a bounded linear map from \(X\) to \(X^*\), with norm at most
\(\norm{D^2f(x)}\). Dividing
\eqref{eq:mvt-first-derivative-remainder} by \(\norm{h}\)
and letting \(h\to0\) proves that it is the \Fre\ derivative
of the map \(Df:B\to X^*\). Thus the second \Fre\ derivative
of \(f\) is \(D^2f\). Its norm-continuity follows from
\eqref{eq:mvt-hessian-lipschitz}, so \(f\in C^2(B)\). Moreover,
\eqref{eq:mvt-hessian-lipschitz} shows that this second derivative
is locally Lipschitz, and therefore
\(f\in C^{2,1}_{\mathrm{loc}}(B)\).

\smallskip
\noindent\emph{Part (ii), Step 4: We will prove that, whenever
\([x,x+h]\subset B\),
\[
 \abs{f(x+h)-f(x)-Df(x)h-\tfrac12D^2f(x)[h^2]}
 \le\tfrac L6\norm{h}^3.
\]}

Let \([x,x+h]\subset B\), and put \(g(s)=f(x+sh)\).
The preceding steps show that \(g\) is an ordinary \(C^2\)
function on an interval containing \([0,1]\), with
\[
 g'(s)=Df(x+sh)h,
 \qquad g''(s)=D^2f(x+sh)[h,h].
\]

For completeness, apply the fundamental theorem of calculus twice:
\[
 g(1)-g(0)-g'(0)
 =\int_0^1\bigl(g'(t)-g'(0)\bigr)\dd t
 =\int_0^1\int_0^t g''(s)\dd s\dd t
 =\int_0^1(1-s)g''(s)\dd s.
\]

The last equality is the change of order in the integral over
the triangle \(0\le s\le t\le1\). It is legitimate because
\(g''\) is continuous. Since \(\int_0^1(1-s)\dd s=1/2\),
subtracting \(g''(0)/2\) gives
\[
 g(1)-g(0)-g'(0)-\tfrac12g''(0)
 =\int_0^1(1-s)\bigl(g''(s)-g''(0)\bigr)\dd s.
\]

In terms of \(f\), this is the exact identity
\begin{equation}\label{eq:mvt-quadratic-remainder-identity}
 \begin{split}
 &f(x+h)-f(x)-Df(x)h-\tfrac12D^2f(x)[h^2]\\
 &\qquad=\int_0^1(1-s)
       (D^2f(x+sh)-D^2f(x))[h,h]\dd s.
 \end{split}
\end{equation}

For \(0\le s\le1\), the integrand satisfies
\[
 |(D^2f(x+sh)-D^2f(x))[h,h]|
 \le\norm{D^2f(x+sh)-D^2f(x)}\norm{h}^2
 \le Ls\norm{h}^3.
\]

Two powers of \(\norm{h}\) come from evaluating a bilinear form
on \((h,h)\); the third comes from its Lipschitz variation along
the segment. Applying this bound to
\eqref{eq:mvt-quadratic-remainder-identity} proves
\[
 \begin{split}
 |f(x+h)-f(x)-Df(x)h-\tfrac12D^2f(x)[h^2]|
 &\le L\norm{h}^3\int_0^1s(1-s)\dd s\\
 &=L\norm{h}^3\left(\frac12-\frac13\right)
 =\frac L6\norm{h}^3.
 \end{split}
\]
This is \eqref{eq:uniform-taylor}, including its precise constant.
The case \(h=0\) is immediate from the same formula.

\smallskip
\noindent\emph{Part (iii): We prove that \(Df\) and \(D^2f\) are bounded on \(B\).}

Let \(y\in B(a,R)\). Since \(\norm{y-a}<R\), comparison with
the center in \eqref{eq:mvt-hessian-lipschitz} gives
\[
 \norm{D^2f(y)}
 \le\norm{D^2f(a)}+L\norm{y-a}
 \le\norm{D^2f(a)}+LR.
\]
Applying \eqref{eq:mvt-gradient-increment} with \(x=a\) and
\(h=y-a\), and then the triangle inequality, yields
\[
 \begin{split}
 \norm{Df(y)}^{*}
 &\le\norm{Df(a)}^{*}+\norm{Df(y)-Df(a)}^{*}\\
 &\le\norm{Df(a)}^{*}+\norm{D^2f(a)}\norm{y-a}
          +\frac L2\norm{y-a}^2\\
 &\le\norm{Df(a)}^{*}+R\norm{D^2f(a)}+\frac L2R^2.
 \end{split}
\]
The right-hand sides are finite and independent of \(y\).
Hence both derivatives are bounded on \(B\).
\end{proof}

The next lemma requires uniform control only at base points at
which the function is still appreciably positive.

\begin{lemma}[Moving-coordinate obstruction]\label{lem:bump}
Let \(2<p<3\). There is no twice \Gat\ differentiable function
\(b:\ell_p\to\R\) with bounded support and \(b(0)=1\) for which
there exist constants \(C,\rho>0\) such that
\begin{equation}\label{eq:superlevel-remainder}
 \abs{b(x+h)-b(x)-Db(x)h-\tfrac12D^2b(x)[h^2]}
 \le C\norm{h}_p^3
\end{equation}
whenever \(b(x)\ge1/2\) and \(\norm{h}_p\le\rho\).
In particular, no nonzero \(G^3\)-smooth bump function on \(\ell_p\)
has a globally bounded third derivative, and no nonzero \(C^2\)
bump function has a globally Lipschitz second derivative.
\end{lemma}
\begin{proof}
Suppose that a function \(b\) and constants \(C,\rho>0\) satisfying
the first part of the statement exist. We will construct points
\(x_j\) whose norms tend to infinity, while their function values
stay strictly above \(1/2\). This will contradict boundedness of
the support.

The construction has two different ingredients. The lengths of
the successive increments are chosen in advance so that their
\(p\)-th powers have divergent sum, but their cubes have a small
finite sum. The coordinates of the increments are then chosen
one at a time, using Lemma~\ref{lem:diagonal} at the current base
point. This makes the contributions of the first and second
derivatives sufficiently small.

\smallskip
\noindent\emph{Step 1: We will prove that one can choose
\(0<t_j\le\rho\) and \(\varepsilon_j>0\) such that
\(\sum_{j=1}^{\infty}t_j^p=\infty\),
\(C\sum_{j=1}^{\infty}t_j^3\le1/8\), and
\(\sum_{j=1}^{\infty}\varepsilon_j=1/4\), so the total permitted
decrease is at most \(3/8\).}

Set
\(S_p=\sum_{j=1}^{\infty}j^{-3/p}\).
Since \(p<3\), we have \(3/p>1\), and the standard convergence
criterion for the numerical series \(\sum_j j^{-q}\), \(q>1\),
gives \(S_p<\infty\). Choose
\[
 \delta=\min\left\{\rho,\,(8CS_p)^{-1/3}\right\}>0,
 \qquad t_j=\delta j^{-1/p}\quad(j\ge1).
\]

Then every step is allowed by the radius in
\eqref{eq:superlevel-remainder}, because
\(0<t_j\le\delta\le\rho\).
Moreover,
\begin{equation}\label{eq:bump-step-sums}
 C\sum_{j=1}^{\infty}t_j^3
 =C\delta^3S_p\le\frac18,
 \qquad
 \sum_{j=1}^{\infty}t_j^p
 =\delta^p\sum_{j=1}^{\infty}\frac1j=\infty.
\end{equation}

The first relation will control the total Taylor error; the
second will force the constructed points to leave every bounded
set. There is no conflict between the two requirements, because
the exponents \(p\) and \(3\) are different and \(p<3\).

We also prescribe an allowance for the linear and quadratic
terms at the \(j\)-th step:
\[
 \varepsilon_j=2^{-j-2},\qquad j\ge1.
\]

The geometric series gives
\(\sum_{j=1}^{\infty}\varepsilon_j
 =\frac14\sum_{j=1}^{\infty}2^{-j}=\frac14\).
 
For \(j\ge0\), let
\(E_j=\sum_{k=1}^{j}\varepsilon_k
       +C\sum_{k=1}^{j}t_k^3\) and \(E_0=0\).
Thus \(E_j\) is an upper bound for the total decrease we will
permit through step \(j\). By \eqref{eq:bump-step-sums},
\begin{equation}\label{eq:bump-total-decrease}
 0\le E_j\le
 \sum_{k=1}^{\infty}\varepsilon_k
       +C\sum_{k=1}^{\infty}t_k^3
 \le\frac14+\frac18=\frac38<\frac12.
\end{equation}

Starting from the value \(b(0)=1\), this leaves a positive margin
above the level \(1/2\) needed to apply the remainder estimate.
The numerical choices \(1/4\) and \(1/8\) are only convenient
allowances; their sum being smaller than \(1/2\) is what matters.

\smallskip
\noindent\emph{Step 2: We will prove that, at each fixed
\(x_{j-1}\), one can choose \(n_j>n_{j-1}\) so that
\(|Db(x_{j-1})h_j|+\tfrac12|D^2b(x_{j-1})[h_j^2]|
\le\varepsilon_j\), where \(h_j=t_je_{n_j}\).}

Suppose that \(x_{j-1}\) and indices \(n_1<\cdots<n_{j-1}\) have
already been chosen. At this step \(x_{j-1}\) is fixed. Twice
\Gat\ differentiability means, in particular, that
\(Db(x_{j-1})\) is a bounded linear functional and
\(D^2b(x_{j-1})\) is a bounded bilinear form on \(\ell_p\).
Applying Lemma~\ref{lem:diagonal} to the normalized disjoint
sequence \((e_n)\) gives
\[
 Db(x_{j-1})e_n\longrightarrow0,
 \qquad
 D^2b(x_{j-1})[e_n,e_n]\longrightarrow0
 \quad(n\to\infty).
\]

The positive numbers \(t_j\) and \(\varepsilon_j\) have already
been fixed. We can therefore choose \(n_j>n_{j-1}\), with
\(n_0=0\) at the first step, so large that simultaneously
\[
 |Db(x_{j-1})e_{n_j}|\le\frac{\varepsilon_j}{2t_j},
 \qquad
 |D^2b(x_{j-1})[e_{n_j},e_{n_j}]|
       \le\frac{\varepsilon_j}{t_j^2}.
\]

To justify the simultaneous choice, each of the two convergences
provides an index beyond which its inequality holds; take \(n_j\)
beyond both indices and beyond \(n_{j-1}\). Multiplication by
\(t_j\) and \(t_j^2/2\), respectively, then gives
\begin{equation}\label{eq:small-coefficients}
 t_j|Db(x_{j-1})e_{n_j}|
 +\frac{t_j^2}{2}|D^2b(x_{j-1})[e_{n_j},e_{n_j}]|
 \le\varepsilon_j=2^{-j-2}.
\end{equation}

Now set
\(h_j=t_je_{n_j}\) and \(x_j=x_{j-1}+h_j\).
The increment has norm \(\norm{h_j}_p=t_j\le\rho\).
The requirement \(n_j>n_{j-1}\) ensures that each new increment
uses a coordinate which has not been used at any earlier step.
No uniform bound for \(Db(x)\) or \(D^2b(x)\) as \(x\) varies is
needed: the coordinate is chosen after the current derivatives
are known.

\smallskip
\noindent\emph{Step 3: We will prove by induction that the
construction is defined for every \(j\ge0\), and that
\(x_j=\sum_{k=1}^{j}t_ke_{n_k}\) satisfies
\(b(x_j)\ge1-E_j\ge5/8>1/2\).}

We prove that the coordinates in Step~2 can be selected
successively and that the resulting points satisfy
\begin{equation}\label{eq:bump-inductive-bound}
 x_j=\sum_{k=1}^{j}t_ke_{n_k},
 \qquad b(x_j)\ge1-E_j\ge\frac58>\frac12
 \quad(j\ge0),
\end{equation}
where the empty sum is \(x_0=0\).
For \(j=0\), this is just \(b(0)=1\) and \(E_0=0\).

Assume that the construction and the estimate have been
established through step \(j-1\). Step~2 supplies a fresh
coordinate \(n_j\) and an increment \(h_j\). At this point the
induction hypothesis already gives
\(b(x_{j-1})\ge5/8>1/2\), and we have
\(\norm{h_j}_p\le\rho\). These are exactly the two hypotheses
for applying \eqref{eq:superlevel-remainder} at the base point
\(x_{j-1}\). Thus
\[
 b(x_j)
 =b(x_{j-1})+t_jDb(x_{j-1})e_{n_j}
 +\frac{t_j^2}{2}
       D^2b(x_{j-1})[e_{n_j},e_{n_j}]+R_j,
 \qquad |R_j|\le Ct_j^3.
\]

Here \(R_j\) denotes the actual difference between the function
value and its displayed quadratic approximation. The estimate
controls its absolute value, so it is valid regardless of the
sign of \(R_j\).

The linear and quadratic terms need not be positive either.
Using their absolute values and \eqref{eq:small-coefficients}
therefore gives the lower bound
\[
 \begin{split}
 b(x_j)
 &\ge b(x_{j-1})
       -t_j|Db(x_{j-1})e_{n_j}|-\frac{t_j^2}{2}
       |D^2b(x_{j-1})[e_{n_j},e_{n_j}]|-Ct_j^3\\
 &\ge b(x_{j-1})-\varepsilon_j-Ct_j^3\\
 &\ge1-E_{j-1}-\varepsilon_j-Ct_j^3
 =1-E_j\ge\frac58.
 \end{split}
\]

The last inequality is \eqref{eq:bump-total-decrease}.
The formula for \(x_j\) in \eqref{eq:bump-inductive-bound}
follows by adding \(t_je_{n_j}\) to the corresponding formula
for \(x_{j-1}\). This completes the induction.

There is no circular use of the remainder hypothesis: the
estimate at \(x_{j-1}\) is justified by the preceding induction
step, and is then used to prove the estimate at \(x_j\).
In particular, we have not assumed \(b(x_j)\ge1/2\) in order
to prove it. Nor do we need the intermediate segment
\([x_{j-1},x_j]\) to lie in the superlevel set. The hypothesis
of this lemma requires only that its base point belong to that
set and that its increment have norm at most \(\rho\).

\smallskip
\noindent\emph{Step 4: We will prove that
\(\norm{x_j}_p^p=\delta^p\sum_{k=1}^{j}1/k\to\infty\),
and hence \(\norm{x_j}_p\to\infty\).}

Since the indices \(n_k\) are all distinct, the coordinates of
the finite sum \(x_j\) do not overlap. By the definition of the
\(\ell_p\)-norm,
\[
 \norm{x_j}_p^p
 =\left\|\sum_{k=1}^{j}t_ke_{n_k}\right\|_p^p
 =\sum_{k=1}^{j}t_k^p
 =\delta^p\sum_{k=1}^{j}\frac1k.
\]

The harmonic partial sums tend to infinity, and \(\delta>0\)
is fixed. Hence \(\norm{x_j}_p\to\infty\). More explicitly,
\[
 \sum_{k=1}^{j}\frac1k\ge\int_1^{j+1}\frac{\dd s}{s}
 =\log(j+1),
 \qquad
 \norm{x_j}_p\ge\delta\bigl(\log(j+1)\bigr)^{1/p}.
\]

Only finite sums have been used to define the points \(x_j\).
We do not form a vector \(\sum_{k=1}^{\infty}t_ke_{n_k}\)
in \(\ell_p\); indeed, that series cannot converge in
\(\ell_p\), because its partial sums have unbounded norms.

\smallskip
\noindent\emph{Step 5: We will prove that bounded support
forces \(b(x_j)=0\) for all sufficiently large \(j\),
contradicting \(b(x_j)\ge5/8\).}

Since \(\supp b\) is bounded, there is a finite \(R_b>0\) such
that
\(\supp b\subset\{x:\norm{x}_p\le R_b\}\).
Every point with \(b(x)\ne0\) belongs to \(\supp b\).
Consequently, from
\(\norm{x}_p>R_b\) it follows that \(b(x)=0\).

Step~4 provides a \(j\) with \(\norm{x_j}_p>R_b\), so
\(b(x_j)=0\). On the other hand, the induction gives
\(b(x_j)\ge5/8\). This contradiction proves the first assertion
of the lemma. Notice that boundedness of the support is
sufficient; compactness of the support is not used.

\smallskip
\noindent\emph{Step 6: We will prove that no nonzero
\(G^3\)-smooth bump \(g:\ell_p\to\R\) can satisfy
\(\sup_{x\in\ell_p}\norm{D^3g(x)}<\infty\).}

Suppose that \(g:\ell_p\to\R\) is a nonzero \(G^3\)-smooth bump
and that
\(\sup_{y\in\ell_p}\norm{D^3g(y)}\le M<\infty\).
The bound \(M\) is taken nonnegative. 

Choose a point \(a\)
with \(\beta=g(a)\ne0\), and define
\(c(x)=\frac{g(a+x)}{\beta}\) for \(x\in\ell_p\).
Then \(c(0)=1\), \(c\) is \(G^3\)-smooth, and
\(\supp c=\supp g-a\) is bounded. Translation of the argument
and division by the fixed nonzero scalar give
\[
 D^kc(x)[h_1,\ldots,h_k]
 =\frac{1}{\beta}D^kg(a+x)[h_1,\ldots,h_k],
 \qquad k=1,2,3.
\]

In particular,
\(\norm{D^3c(x)}\le\frac{M}{|\beta|}\) for \(x\in\ell_p\).
The absolute value in the denominator handles either sign of
\(\beta\); nonnegativity of the original bump is not required.

For arbitrary \(x,h\in\ell_p\), apply
Lemma~\ref{lem:mvt} to \(c\) on an open ball containing the
segment \([x,x+h]\). The bound \(M/|\beta|\) is independent
of that ball, so \eqref{eq:uniform-taylor} gives
\begin{equation}\label{eq:bump-normalized-remainder}
 |c(x+h)-c(x)-Dc(x)h-\tfrac12D^2c(x)[h^2]|
 \le\frac{M}{6|\beta|}\norm{h}_p^3
 \qquad(x,h\in\ell_p).
\end{equation}

Thus the first part of the lemma applies to \(c\), for example
with
\(\rho=1\) and \(C=\max\left\{1,\frac{M}{6|\beta|}\right\}>0\).
The estimate holds at every base point, and in particular
whenever \(c(x)\ge1/2\). We have obtained the function excluded
by Steps~1--5, a contradiction. This proves the assertion about
globally bounded third derivatives.

\smallskip
\noindent\emph{Step 7: We will prove that no nonzero
\(C^2\) bump \(g:\ell_p\to\R\) admits a constant \(M\ge0\)
such that \(\norm{D^2g(y)-D^2g(z)}\le M\norm{y-z}_p\)
for all \(y,z\in\ell_p\).}

Suppose instead that \(g\) is a nonzero \(C^2\) bump and that
some \(M\ge0\) satisfies
\[
 \norm{D^2g(y)-D^2g(z)}
 \le M\norm{y-z}_p
 \qquad(y,z\in\ell_p).
\]
Use the same normalization \(c(x)=g(a+x)/\beta\), where
\(\beta=g(a)\ne0\). It preserves \(C^2\)-regularity and
bounded support, and yields
\[
 c(0)=1,\qquad
 \norm{D^2c(x)-D^2c(y)}
 \le\frac{M}{|\beta|}\norm{x-y}_p
 \qquad(x,y\in\ell_p).
\]
We do not assume that \(D^3g\), or \(D^3c\), exists.
The integral identity
\eqref{eq:mvt-quadratic-remainder-identity} needs only
\(C^2\)-regularity, as its derivation in part~(ii), Step~4,
of the proof of Lemma~\ref{lem:mvt} shows. Applying that identity to \(c\)
and using its Lipschitz second derivative gives
\[
 \begin{split}
 &|c(x+h)-c(x)-Dc(x)h-\tfrac12D^2c(x)[h^2]|\\
 &\qquad\le\int_0^1(1-s)
       \norm{D^2c(x+sh)-D^2c(x)}\norm{h}_p^2\dd s\\
 &\qquad\le\frac{M}{|\beta|}\norm{h}_p^3
       \int_0^1s(1-s)\dd s
 =\frac{M}{6|\beta|}\norm{h}_p^3.
 \end{split}
\]
We therefore recover \eqref{eq:bump-normalized-remainder}
without using any third derivative. The first part of the lemma,
with the same positive choice of \(C\) and \(\rho=1\), again
gives a contradiction. This proves the final assertion.
\end{proof}

\begin{lemma}[Localization at a strongly exposed point]\label{lem:cap}
Let \(X\) be a real Banach space and \(\tnorm{\cdot}\) an equivalent norm.
Suppose that \(x_0\in S_{\tnorm{\cdot}}\) is strongly exposed and that \(\tnorm{\cdot}\) is
\(G^3\)-smooth on a ball about \(x_0\), with its third derivative
bounded there. Then \(X\) admits a nonzero \(G^3\)-smooth bump
function with globally bounded third derivative.
\end{lemma}
\begin{proof}
Write \(\norm{\cdot}\) for the original norm of \(X\), which will
be used for distances, open balls, and operator norms of
derivatives. The unit ball and sphere of the equivalent norm
\(\tnorm{\cdot}\) are denoted by \(B_{\tnorm{\cdot}}\) and \(S_{\tnorm{\cdot}}\), respectively.
The construction will use a product
\(b(x)=\chi(\tnorm{x})\psi(f(x))\): one scalar factor restricts the
size of \(\tnorm{x}\), and the other selects points near the
supporting hyperplane at \(x_0\). Strong exposure will force
their simultaneous support into the prescribed neighborhood of
\(x_0\).

\smallskip
\noindent\emph{Step 1: We will prove that one can choose
\(0<r<\norm{x_0}/2\), \(L\ge0\), and a functional \(f\)
strongly exposing \(x_0\), so that
\(U=B(x_0,r)\subset X\setminus\{0\}\),
\(\tnorm{\cdot}\in G^3(U)\), \(\sup_{x\in U}\norm{D^3\tnorm{\cdot}(x)}\le L\),
and \(\tnorm{f}^{*}=f(x_0)=1\).}

By hypothesis, after shrinking the given neighborhood if
necessary, we can choose \(r>0\) and \(L\ge0\) such that
\(U=B(x_0,r)\subset X\setminus\{0\}\), \(\tnorm{\cdot}\in G^3(U)\), and \(\sup_{x\in U}\norm{D^3\tnorm{\cdot}(x)}\le L\).
We may additionally require \(r<\norm{x_0}/2\), since
\(\tnorm{x_0}=1\) implies \(x_0\ne0\). Shrinking the ball does not
affect the differentiability assumption or invalidate its
derivative bound.

Choose a functional \(f\in X^*\) strongly exposing \(x_0\),
normalized so that
\(\tnorm{f}^{*}=1\) and \(f(x_0)=1\).
Thus \(|f(x)|\le \tnorm{x}\) for every \(x\in X\), and in
particular \(f(y)\le1\) for \(y\in B_{\tnorm{\cdot}}\).

Strong exposure means that
if \(y_j\in B_{\tnorm{\cdot}}\) and \(f(y_j)\to1\), then \(\norm{y_j-x_0}\to0\).

The convergence may be expressed in the original norm because
that norm and \(\tnorm{\cdot}\) are equivalent. We will also use the
original dual norm
\(F_0=\norm{f}^{*}
 =\sup_{\norm{h}\le1}|f(h)|<\infty\).
It need not equal \(\tnorm{f}^{*}=1\); the derivative bounds
below will keep this distinction explicit.

\smallskip
\noindent\emph{Step 2: We will prove that there exists
\(\eta\in(0,1)\) such that
\(\{y\in B_{\tnorm{\cdot}}:f(y)>1-\eta\}\subset B(x_0,r/4)\).}

There exists \(0<\eta<1\) such that
\begin{equation}\label{eq:cap-small-slice}
 \tnorm{y}\le1,\quad f(y)>1-\eta
 \quad\Longrightarrow\quad\norm{y-x_0}<r/4.
\end{equation}

Here is the sequential justification. If no such \(\eta\)
existed, for each \(j\ge2\) we could choose \(y_j\in B_{\tnorm{\cdot}}\)
with
\(f(y_j)>1-\frac1j\) and \(\norm{y_j-x_0}\ge r/4\).
Since also \(f(y_j)\le1\), we would have \(f(y_j)\to1\).
Strong exposure would imply \(\norm{y_j-x_0}\to0\), contrary
to the lower bound \(r/4\). This proves
\eqref{eq:cap-small-slice}; reducing \(\eta\), if necessary,
ensures \(0<\eta<1\).

Thus a sufficiently thin slice of the unit ball, cut out by
the single functional \(f\), is contained in a small norm ball
about \(x_0\). The radius \(r/4\) leaves room for a subsequent
small change in the radial coordinate.

\smallskip
\noindent\emph{Step 3: We will prove that one can choose
\(\varepsilon>0\) and \(\chi,\psi\in C_c^\infty(\R)\), with
\(\chi(1)=\psi(1)=1\) and
\(\supp\chi,\supp\psi\subset(1-\varepsilon,1+\varepsilon)\),
so that \(b(x)=\chi(\tnorm{x})\psi(f(x))\) is continuous on \(X\)
and satisfies \(b(x_0)=1\).}

Choose a positive \(\varepsilon\) satisfying
\(0<\varepsilon<
 \min\left\{\frac12,\frac{\eta}{4},
             \frac{r}{8\norm{x_0}+2r}\right\}\).
The first two restrictions give
\begin{equation}\label{eq:cap-ratio}
 \frac{1-\varepsilon}{1+\varepsilon}
 =1-\frac{2\varepsilon}{1+\varepsilon}
 >1-2\varepsilon>1-\eta.
\end{equation}

The last restriction gives a quantitative margin inside
\(B(x_0,r/2)\). 

Namely, if
\(r_1=(1+\varepsilon)\frac r4+\varepsilon\norm{x_0}\),
then
\begin{equation}\label{eq:cap-margin}
 r_1=\frac r4+\varepsilon\left(\frac r4+\norm{x_0}\right)
 <\frac r4+\frac r8=\frac{3r}{8}<\frac r2.
\end{equation}

All constants have now been fixed, before choosing the
function \(b\).

Take \(\chi,\psi\in C_c^\infty(\R)\) such that
\(0\le\chi,\psi\le1\), \(\chi(1)=\psi(1)=1\), and \(\supp\chi,\supp\psi\subset(1-\varepsilon,1+\varepsilon)\).

For completeness, an explicit choice is obtained from
\[
 \theta(s)=
 \begin{cases}
  \exp\!\left(1-\dfrac{1}{1-s^2}\right),&|s|<1,\\
  0,&|s|\ge1,
 \end{cases}
 \qquad
 \chi(s)=\psi(s)=\theta\!\left(\frac{2(s-1)}{\varepsilon}\right).
\]

The function \(\theta\) is \(C^\infty\) on \(\R\):
all derivatives from inside \((-1,1)\) tend to zero at the
endpoints, since the product of the exponential factor with
any power of \((1-s^2)^{-1}\) tends to zero. They therefore
match the derivatives of the zero function outside this interval.
Also \(\theta(0)=1\), \(0\le\theta\le1\), and
\(\supp\chi=\supp\psi
 =[1-\varepsilon/2,\,1+\varepsilon/2]
 \subset(1-\varepsilon,1+\varepsilon)\).
It is harmless to choose the two cutoffs equal; their different
names distinguish their roles in the composition.

Define on the whole space \(X\)
\(b(x)=\chi(\tnorm{x})\psi(f(x))\).
The norm \(\tnorm{\cdot}\) and the functional \(f\) are continuous on \(X\),
so \(b\) is continuous. Moreover,
\(b(x_0)=\chi(1)\psi(1)=1\).
Thus \(b\) is nonzero. We next locate its entire closed support.

\smallskip
\noindent\emph{Step 4: We will prove that
\(\supp b\subset B(x_0,r_1)\subset B(x_0,r/2)\subset U\),
where \(r_1<r/2\) is the radius fixed in Step~3. In particular,
\(b\) has bounded support.}

Consider the set
\(K=\{x\in X:\tnorm{x}\in\supp\chi,\ f(x)\in\supp\psi\}\).
Both scalar supports are closed and both maps \(\tnorm{\cdot},f\) are
continuous. Hence \(K\) is closed in \(X\).
If \(b(x)\ne0\), both factors in its definition are nonzero,
so \(x\in K\). Taking the closure of this inclusion gives
\[
 \supp b=\overline{\{x:b(x)\ne0\}}\subset K.
\]

This closed-set argument ensures that boundary points of the
support are included in the localization, even if \(b\) itself
vanishes at those points.

Let \(x\in K\). The support restrictions give
\[
 1-\varepsilon\le \tnorm{x}\le1+\varepsilon,\qquad
 1-\varepsilon\le f(x)\le1+\varepsilon.
\]

In particular, \(\tnorm{x}>0\). We may therefore normalize \(x\)
by setting
\(y=\frac{x}{\tnorm{x}}\).
Homogeneity of \(\tnorm{\cdot}\) and linearity of \(f\) give
\(\tnorm{y}=1\) and \(f(y)=\frac{f(x)}{\tnorm{x}}
 \ge\frac{1-\varepsilon}{1+\varepsilon}>1-\eta\),
where the last inequality is \eqref{eq:cap-ratio}. By
\eqref{eq:cap-small-slice}, \(\norm{y-x_0}<r/4\).
To pass back from \(y\) to \(x\), use the exact identity
\[
 x-x_0=\tnorm{x}(y-x_0)+(\tnorm{x}-1)x_0.
\]

It follows that
\[
 \norm{x-x_0}
 \le \tnorm{x}\norm{y-x_0}+|\tnorm{x}-1|\norm{x_0}
 <(1+\varepsilon)\frac r4+\varepsilon\norm{x_0}=r_1.
\]

Consequently,
\begin{equation}\label{eq:cap-support-localization}
 \supp b\subset K\subset B(x_0,r_1)
 \subset B(x_0,r/2)\subset U.
\end{equation}

In particular, the support is bounded, so \(b\) is a bump
function. The strict inequality \(r_1<r/2\) in
\eqref{eq:cap-margin} leaves a positive distance between the
support and the complement of \(B(x_0,r/2)\). This margin will
also handle points on the boundary of the region where we
estimate the derivatives.

\smallskip
\noindent\emph{Step 5: We will prove that
\(b(x)=\chi(\tnorm{x})\psi(f(x))\) belongs to \(G^3(U)\),
and derive its third derivative by the chain and product rules
in the operator-norm sense.}

On \(U\), put
\(g(x)=\chi(\tnorm{x})\) and \(q(x)=\psi(f(x))\), so that \(b(x)=g(x)q(x)\).
The function \(q\) is \(C^\infty\) on all of \(X\), because
it is the composition of a smooth scalar function with a
continuous linear functional. The function \(g\) is \(G^3\)
on \(U\). To justify this last claim in the recursive
operator-norm convention, first apply Lemma~\ref{lem:mvt}
to \(\tnorm{\cdot}\). It gives \(\tnorm{\cdot}\in C^2(U)\). The ordinary chain rule
therefore gives the first two derivative identities below.
Differentiating the second identity along a fixed direction
in operator norm is then legitimate: \(D^2\tnorm{\cdot}\) has its
directional derivative in operator norm by the \(G^3\)
assumption, and the scalar factors and \(D\tnorm{\cdot}\) are
norm-continuous and differentiable as required.
Multiplication of a scalar and a bounded form, and products
of bounded linear functionals, are continuous operations in
the corresponding operator norms.

For \(x\in U\) and \(h,k,l\in X\), the resulting formulas are
\(Dg(x)h=\chi'(\tnorm{x})D\tnorm{\cdot}(x)h\)
and
\[
 D^2g(x)[h,k]
 =\chi''(\tnorm{x})D\tnorm{\cdot}(x)h\,D\tnorm{\cdot}(x)k
       +\chi'(\tnorm{x})D^2\tnorm{\cdot}(x)[h,k].
\]

The third derivative is
\begin{equation}\label{eq:cap-chain-third}
 \begin{split}
 D^3g(x)[h,k,l]
 &=\chi'''(\tnorm{x})D\tnorm{\cdot}(x)h\,D\tnorm{\cdot}(x)k\,D\tnorm{\cdot}(x)l\\
 &\quad+\chi''(\tnorm{x})\bigl(
       D^2\tnorm{\cdot}(x)[h,k]\,D\tnorm{\cdot}(x)l\\
 &\quad+D^2\tnorm{\cdot}(x)[h,l]\,D\tnorm{\cdot}(x)k
 +D^2\tnorm{\cdot}(x)[k,l]\,D\tnorm{\cdot}(x)h\bigr)\\
 &\quad+\chi'(\tnorm{x})D^3\tnorm{\cdot}(x)[h,k,l].
 \end{split}
\end{equation}

For the second factor, linearity of \(f\) eliminates all
higher derivatives of the inner map, and gives
\begin{equation}\label{eq:cap-linear-composition}
 D^jq(x)[h_1,\ldots,h_j]
 =\psi^{(j)}(f(x))\prod_{i=1}^{j}f(h_i),
 \qquad j=1,2,3.
\end{equation}

Thus \(b=gq\) is \(G^3\) on \(U\). For clarity, its full
third-derivative product rule is
\begin{equation}\label{eq:cap-product-third}
 \begin{split}
 D^3b(x)[h,k,l]
 &=q(x)D^3g(x)[h,k,l]\\
 &\quad+D^2g(x)[h,k]\,Dq(x)l
 +D^2g(x)[h,l]\,Dq(x)k
 +D^2g(x)[k,l]\,Dq(x)h\\
 &\quad+Dg(x)h\,D^2q(x)[k,l]
 +Dg(x)k\,D^2q(x)[h,l]
 +Dg(x)l\,D^2q(x)[h,k]\\
 &\quad+g(x)D^3q(x)[h,k,l].
 \end{split}
\end{equation}

The three terms of each mixed type account for the three
ways of assigning one of the directions to one factor and
the other two directions to the other factor. These formulas
use existence of the third \Gat\ derivative in the stated
operator-norm sense; they do not require continuity of
\(D^3\tnorm{\cdot}\) as a function of \(x\).

\smallskip
\noindent\emph{Step 6: We will prove that there is a finite
constant \(K_0\), independent of \(x\in B(x_0,r/2)\),
such that \(\norm{D^3b(x)}\le K_0\).}

Let \(V=B(x_0,r/2)\). The explicit lower-derivative bounds
from Lemma~\ref{lem:mvt} give, for \(x\in V\),
\(\norm{D\tnorm{\cdot}(x)}^{*}\le A_1\), \(\norm{D^2\tnorm{\cdot}(x)}\le A_2\), and \(\norm{D^3\tnorm{\cdot}(x)}\le L\),
where we may take the finite constants
\(A_1=\norm{D\tnorm{\cdot}(x_0)}^{*}
          +\frac r2\norm{D^2\tnorm{\cdot}(x_0)}+\frac{Lr^2}{8}\) and
\(A_2=\norm{D^2\tnorm{\cdot}(x_0)}+\frac{Lr}{2}\).

For \(j=0,1,2,3\), set
\(c_j=\sup_{s\in\R}|\chi^{(j)}(s)|\) and \(d_j=\sup_{s\in\R}|\psi^{(j)}(s)|\),
where the superscript \(0\) denotes the function itself.
These constants are finite because the scalar cutoffs are
smooth and compactly supported. For the choices above,
\(c_0=d_0=1\). From the chain formulas, for every \(x\in V\),
\[
 \begin{split}
 |g(x)|&\le c_0,\\
 \norm{Dg(x)}^{*}&\le c_1A_1,\\
 \norm{D^2g(x)}&\le c_2A_1^2+c_1A_2,\\
 \norm{D^3g(x)}&\le c_3A_1^3+3c_2A_1A_2+c_1L.
 \end{split}
\]

For example, each of the three terms multiplied by
\(\chi''(\tnorm{x})\) in \eqref{eq:cap-chain-third} is bounded
in modulus by \(c_2A_2A_1\) on unit arguments. This explains
the factor \(3\) in the last line. Similarly,
\eqref{eq:cap-linear-composition} gives
\[
 |q(x)|\le d_0,\qquad
 \norm{Dq(x)}^{*}\le d_1F_0,\qquad
 \norm{D^2q(x)}\le d_2F_0^2,\qquad
 \norm{D^3q(x)}\le d_3F_0^3.
\]

Taking absolute values in \eqref{eq:cap-product-third},
and then taking the supremum over unit \(h,k,l\), yields
\[
 \begin{split}
 \norm{D^3b(x)}
 &\le |q(x)|\norm{D^3g(x)}
       +3\norm{D^2g(x)}\norm{Dq(x)}^{*}\\
 &\quad+3\norm{Dg(x)}^{*}\norm{D^2q(x)}
       +|g(x)|\norm{D^3q(x)}.
 \end{split}
\]

Consequently, a single bound valid at every \(x\in V\) is
\begin{equation}\label{eq:cap-global-constant}
 \begin{split}
 \norm{D^3b(x)}\le K_0
 &=d_0(c_3A_1^3+3c_2A_1A_2+c_1L)
 +3d_1F_0(c_2A_1^2+c_1A_2)\\
 &\quad+3d_2F_0^2c_1A_1+d_3F_0^3c_0<\infty.
 \end{split}
\end{equation}

The cutoffs and their derivative bounds may depend on the
chosen slice and neighborhood. This causes no problem: they
have been fixed once, so \(K_0\) is independent of the base
point and of the directions in the derivative.

\smallskip
\noindent\emph{Step 7: We will prove that \(b\in G^3(X)\)
and \(\sup_{x\in X}\norm{D^3b(x)}\le K_0\), so that
\(b\) is the required nonzero bump with globally bounded
third derivative.}

We have proved \(G^3\)-smoothness on \(U\) and a common
third-derivative bound on \(V\). Let \(z\in X\setminus V\).
By \eqref{eq:cap-support-localization}, \(z\notin\supp b\).
Since \(\supp b\) is closed, its complement is open; hence
\(b\) vanishes identically on a neighborhood of \(z\).
On that neighborhood all derivatives of \(b\) exist and are
zero, irrespective of the regularity of \(\tnorm{\cdot}\) there.
This argument includes points on the boundary of \(V\).

For completeness, the positive margin in
\eqref{eq:cap-margin} also gives a direct neighborhood:
with \(d=r/2-r_1>0\), if \(z\notin V\) and
\(\norm{w-z}<d/2\), then
\[
 \norm{w-x_0}
 \ge\norm{z-x_0}-\norm{w-z}
 >\frac r2-\frac d2>r_1.
\]
Thus \(w\notin\supp b\), and \(b=0\) on \(B(z,d/2)\).
At a boundary point of \(\supp b\) itself, we instead use
the local chain and product formulas, since that point lies
inside \(V\subset U\). These two observations cover all
potential boundary issues without differentiating \(\tnorm{\cdot}\)
outside its given region of \(G^3\)-smoothness.

It follows that \(b\in G^3(X)\), and
\(\sup_{x\in X}\norm{D^3b(x)}
 \le K_0\),
because \eqref{eq:cap-global-constant} applies on \(V\) and
the third derivative is zero on \(X\setminus V\).
Together with continuity, bounded support, and \(b(x_0)=1\),
this proves the lemma.
\end{proof}

De~Bernardi, Preti and Somaglia use a weighted sum of squares of
functionals to obtain strict convexity while retaining first-order
\Gat\ smoothness \cite[proof of Theorem~3.3]{DeBernardiPretiSomaglia}.
Here the perturbation must also preserve the recursive third derivative.

\begin{lemma}[Compact Hilbertian perturbation]\label{lem:kadec}
Let \(1<p<\infty\), let \(\tnorm{\cdot}\) be an equivalent Kadec norm on
\(\ell_p\), and put
\[
 H(x)=\left(\sum_{k=1}^{\infty}2^{-k}x_k^2\right)^{1/2},
 \qquad \tnorm{x}_H=\tnorm{x}+H(x).
\]
Then 
\begin{enumerate}
	\item \(\tnorm{\cdot}_H\) is an equivalent strictly convex Kadec norm.
	\item Every point of its unit sphere is strongly exposed by any norming
	functional.
	\item \(H\) is \(C^\infty\) on
	\(\ell_p\setminus\{0\}\).
\end{enumerate}
Hence \(G^3\)-smoothness of \(\tnorm{\cdot}\), when
present, is preserved.
\end{lemma}

\begin{proof}
The map \(T:\ell_p\to\ell_2\),
\(Tx=(2^{-k/2}x_k)_k\), is bounded, injective and compact. Indeed,
its finite-coordinate truncations converge in operator norm, since
\[
 \sum_{k>r}2^{-k}|x_k|^2
 \le\norm{x}_p^2\sum_{k>r}2^{-k}.
\]
Thus \(H=\norm{T(\cdot)}_2\) is strictly convex in the sense that
equality in its triangle inequality for two nonzero vectors forces
them to be positively proportional. Adding \(H\) makes \(\tnorm{\cdot}_H\)
strictly convex, while equivalence follows from boundedness of \(T\).
Since \(Q(x)=\norm{Tx}_2^2>0\) away from zero,
\(H=\sqrt Q\) is \(C^\infty\) there.

Suppose \(x_n\rightharpoonup x\) and
\(\tnorm{x_n}_H=\tnorm{x}_H=1\).
Compactness gives \(H(x_n)\to H(x)\), and hence
\(\tnorm{x_n}\to \tnorm{x}>0\). The normalized vectors
\(x_n/\tnorm{x_n}\) converge weakly to \(x/\tnorm{x}\). The Kadec property of
\(\tnorm{\cdot}\) now gives norm convergence. This proves the Kadec property
of \(\tnorm{\cdot}_H\). Sequential verification suffices here, since
bounded weak subsets of \(\ell_p\) are metrizable.

Fix \(x\in S_{\tnorm{\cdot}_H}\) and a norming functional \(f\),
so that \(\tnorm{f}_H^{*}=f(x)=1\).
If \(y_n\in B_{\tnorm{\cdot}_H}\) and \(f(y_n)\to1\), reflexivity
provides a weakly convergent subsequence of every subsequence.
Every weak cluster point
belongs to the face
\(\{y:\tnorm{y}_H\le1,\ f(y)=1\}\), which is the singleton
\(\{x\}\) by strict convexity. Thus \(y_n\rightharpoonup x\), and
\(\tnorm{y_n}_H\to1\). Normalization and the Kadec property
give \(\norm{y_n-x}\to0\), proving strong exposure.
\end{proof}

\section{Proofs of the main results}\label{sec:proofs}

\begin{proof}[Proof of Theorem~\ref{thm:baire}]
1.) 
The case \(X=\{0\}\) is immediate, so assume that \(X\ne\{0\}\).
Openness of \(G_f(U)\) follows from its definition. To prove density,
start with an arbitrary nonempty open subset \(O\) of \(U\).
Choose a closed ball \(K\) of positive radius contained in \(O\)
and with positive distance from \(X\setminus U\).
Let \(H\) be a closed direction ball of positive finite radius.
For sufficiently large \(m_0\), all points \(x+jth\),
\(j=0,1,2,3\), belong to \(U\) whenever \(x\in K\), \(h\in H\)
and \(0<t<1/m_0\).

Define the third difference by
\[
 \Delta_{th}^3 f(x)=f(x+3th)-3f(x+2th)+3f(x+th)-f(x).
\]
Taylor's formula along the fixed line through \(x\) in direction
\(h\) gives
\begin{equation}\label{eq:finite-diff-limit}
 \frac{\Delta_{th}^3 f(x)}{t^3}
 \longrightarrow D^3f(x)[h^3].
\end{equation}
For integers \(m\ge m_0\), put
\[
 E_m=\{(x,h)\in K\times H:
        |\Delta_{th}^3 f(x)|\le mt^3
        \text{ for every rational }t\in(0,1/m)\}.
\]
Continuity of \(f\) makes each \(E_m\) closed in \(K\times H\).
For each fixed pair \((x,h)\), the limit in
\eqref{eq:finite-diff-limit} bounds the quotient for all sufficiently
small \(t\); increasing \(m\) controls both the bound and the
admissible interval. Thus the sets \(E_m\) cover the complete
metric space \(K\times H\).
By Baire's theorem, some \(E_m\) has nonempty relative interior.
Since the interiors of \(K\) and \(H\) are dense in their
respective closed balls, there are ambient open balls \(V\subset O\)
and \(W=B(h_0,r)\) with \(V\times W\subset E_m\).
Taking the limit through rational \(t\) gives
\begin{equation}\label{eq:cubic-ball}
 |D^3f(x)[h^3]|\le m\qquad(x\in V,\ h\in W).
\end{equation}

Fix \(x\in V\), a unit vector \(u\), and set
\(P_x(h)=D^3f(x)[h^3]\) and \(a=r/4\).
The four points \(h_0+ja u\), \(j=0,1,2,3\), lie in \(W\).
The third difference of this homogeneous cubic polynomial is
\[
 P_x(h_0+3au)-3P_x(h_0+2au)
 +3P_x(h_0+au)-P_x(h_0)=6a^3P_x(u).
\]
Consequently \eqref{eq:cubic-ball} yields
\(|P_x(u)|\le8m/(6a^3)\). Polarization of a real symmetric
trilinear form then gives
\(\norm{D^3f(x)}\le(3^3/3!)\sup_{\norm{u}\le1}|P_x(u)|
\le6m/a^3\); see the standard polarization formula in \cite{HJ}.
The bound is independent of \(x\in V\), so \(V\subset G_f(U)\).
As \(O\) was arbitrary, \(G_f(U)\) is dense in \(U\).

2.)
Fix \(x\in G_f(U)\). Choose an open ball \(B\) about \(x\) contained
in \(U\) such that \(\sup_{y\in B}\norm{D^3f(y)}<\infty\). Each
point of \(B\) has a neighborhood in \(B\), so \(B\subset G_f(U)\).
Lemma~\ref{lem:mvt} gives \(f\in C^{2,1}_{\mathrm{loc}}(B)\).
Such balls cover \(G_f(U)\), proving the assertion.

3.)
Let \(B\subset U\) be an open ball on which \(\norm{D^3f}\le L\).
For every segment \([x,x+h]\subset B\), Lemma~\ref{lem:mvt}
applied to \(f|_B\) bounds the quadratic Taylor remainder by
\(L\norm{h}^3/6\), as claimed in \eqref{eq:uniform-taylor}.
\end{proof}

\begin{proof}[Proof of Theorem~\ref{thm:kadec}]
1.)
Theorem~\ref{thm:baire}\textup{(\ref{it:baire-open})} makes \(G_{\tnorm{\cdot}}\) open and dense
in \(\ell_p\setminus\{0\}\). Homogeneity gives
\(D^3\tnorm{\cdot}(\lambda x)=\lambda^{-2}D^3\tnorm{\cdot}(x)\)
for \(\lambda>0\). Scaling a neighborhood on which the third
derivative is bounded therefore shows that
\(G_{\tnorm{\cdot}}\) is invariant under positive scaling.
Hence \(G_{\tnorm{\cdot}}\cap S_{\tnorm{\cdot}}\) is relatively
open and dense in \(S_{\tnorm{\cdot}}\): for density, approximate
a point of the sphere by points of \(G_{\tnorm{\cdot}}\) and
apply the continuous normalization \(x\mapsto x/\tnorm{x}\).

If \(x_0\in G_{\tnorm{\cdot}}\cap S_{\tnorm{\cdot}}\) were
strongly exposed, a sufficiently small ball around \(x_0\)
would satisfy the hypotheses of Lemma~\ref{lem:cap}. This would
produce a nonzero \(G^3\)-smooth bump with globally bounded
third derivative, contrary to Lemma~\ref{lem:bump}.

2.)
By assertion~\ref{it:kadec-regular}, every strongly exposed point of
the unit ball belongs to \(S_{\tnorm{\cdot}}\setminus G_{\tnorm{\cdot}}\).
This set is relatively closed and has empty relative interior because
\(G_{\tnorm{\cdot}}\cap S_{\tnorm{\cdot}}\) is relatively open and dense
in the sphere.

3.)
Suppose that \(\tnorm{\cdot}\) were Kadec.
Lemma~\ref{lem:kadec} gives the equivalent \(G^3\)-smooth norm
\(\tnorm{x}_H=\tnorm{x}+H(x)\), every point of whose unit sphere
is strongly exposed. Applying assertion~\ref{it:kadec-regular} to this new norm contradicts
the nonemptiness of its relatively open dense regular set.
Therefore \(\tnorm{\cdot}\) is not Kadec.
\end{proof}

\begin{proof}[Proof of Theorem~\ref{thm:counter}]
1.)
The map \((a,b)\mapsto(a^m+b^m)^{1/m}\) is convex and
nondecreasing in each coordinate on the nonnegative quadrant.
Since \(Q\) and \(M\) are nonnegative and convex, their
composition in \eqref{eq:F} is convex. It is even and continuous,
and satisfies
\begin{equation}\label{eq:F-bounds}
 \max\{Q(x),M(x)\}\le \Psi(x)\le Q(x)+M(x).
\end{equation}
Both \(Q\) and \(M\) are \(C^2\), and both are positive off
zero, so \(\Psi\) is \(C^2\) there. At zero,
\eqref{eq:F-bounds} gives \(\Psi(x)=O(\norm{x}_p^2)\), whence
\(D\Psi(0)=0\) in the \Fre\ sense. The partial derivatives of
\((a^m+b^m)^{1/m}\) belong to \([0,1]\) for positive \(a,b\).
It follows that
\[
 \norm{D\Psi(x)}^{*}\le\norm{DQ(x)}^{*}+\norm{DM(x)}^{*}
 \le2\norm{x}_p+p\norm{x}_p^{p-1}\longrightarrow0,
\]
which also proves \(C^1\)-regularity at zero.

2.)
For \(x\ne0\), write \(a=Q(x)\),
\(b=M(x)\), \(A=(a/\Psi(x))^{m-1}\), and
\(B=(b/\Psi(x))^{m-1}\). Then \(D\Psi(x)=A\,DQ(x)+B\,DM(x)\),
where \(0\le B\le1\) and
\(0\le1-A\le\min\{1,(b/a)^m\}\).
Indeed, the last inequality follows by applying the scalar mean
value theorem to \(s\mapsto(1+s)^{-(m-1)/m}\) for \(s\ge0\).
Writing \(Q(x)=\norm{Tx}_2^2\), with \(T\) as in
Lemma~\ref{lem:kadec}, gives
\(\norm{DQ(x)}^{*}\le2\sqrt{Q(x)}\), since \(\norm{T}\le1\).
If \(a\le b\), then \(\sqrt a\le\sqrt b\); if \(a\ge b\),
then \(\sqrt a(b/a)^m=\sqrt b(b/a)^{m-1/2}\le\sqrt b\).
Consequently,
\begin{equation}\label{eq:F-frechet-second}
 \norm{D\Psi(x)-DQ(x)}^{*}
 \le2\sqrt{M(x)}+p\norm{x}_p^{p-1}
 =2\norm{x}_p^{p/2}+p\norm{x}_p^{p-1}
 =o(\norm{x}_p).
\end{equation}
Because \(DQ(x)=D^2Q[x,\cdot]\) and both exponents on the
right exceed one, this proves that \(D\Psi\) is \Fre\
differentiable at zero, with \(D^2\Psi(0)=D^2Q\).
Together with \(C^2\)-regularity off zero, it proves twice
\Fre\ differentiability everywhere.

Put \(\alpha=m(p-2)>1\). For \(t>0\), set
\(\lambda=t^{p-2}\) and
\(G_\lambda(y)=(Q(y)^m+\lambda^mM(y)^m)^{1/m}\).
Homogeneity yields
\begin{equation}\label{eq:F-scaling}
 \Psi(ty)=t^2G_\lambda(y).
\end{equation}
Fix \(h\ne0\). On a sufficiently small ball about \(h\),
\(Q\) is bounded below by a positive constant, and \(Q,M\)
and their derivatives through order two are bounded.
Consequently \(R=(M/Q)^m\) and its first two derivatives are
bounded there. Write \(G_\lambda-Q=Q\theta(\lambda^mR)\),
where \(\theta(s)=(1+s)^{1/m}-1\). Since \(\theta(0)=0\)
and \(\theta',\theta''\) are bounded near zero, the function
\(\theta(\lambda^mR)\), its first derivative and its second
derivative are all \(O_h(\lambda^m)\) on that ball.
Indeed, every term obtained by differentiating the composition
contains at least one factor \(\lambda^m\); the term with two
first derivatives of \(R\) contains \(\lambda^{2m}\).
The product rule therefore gives
\begin{equation}\label{eq:C2-scaling}
 \norm{D^jG_\lambda(h)-D^jQ(h)}
 =O_h(\lambda^m),\qquad j=0,1,2.
\end{equation}
In this formula and the next one, the norm denotes absolute
value for \(j=0\), the original dual norm \(\norm{\cdot}^{*}\)
for \(j=1\), and the bilinear operator norm for \(j=2\).
The constants may depend on \(h\); no bound uniform over
directions is claimed.

Differentiating \eqref{eq:F-scaling} with respect to \(y\),
and using the quadratic homogeneity of \(Q\), yields
\begin{equation}\label{eq:F-derivative-errors}
 \norm{D^j(\Psi-Q)(th)}
 =O_h(t^{2+\alpha-j}),\qquad j=0,1,2.
\end{equation}
We have already proved \(D^2\Psi(0)=D^2Q\), so the second
derivative is defined throughout a neighborhood of zero.
The case \(j=2\) of \eqref{eq:F-derivative-errors} gives
\[
 \frac{\norm{D^2\Psi(th)-D^2\Psi(0)}}{t}
 =O_h(t^{\alpha-1})\longrightarrow0.
\]
Evenness of \(\Psi\) makes \(D\Psi\) odd and \(D^2\Psi\) even away
from zero, so the same limits hold for negative \(t\).
The direction \(h=0\) is immediate. This proves the recursive
operator-norm third derivative in \eqref{eq:jets-F}.

3.)
By \eqref{eq:F-bounds}, \(\Psi(h)\ge M(h)=\norm{h}_p^p\) for every
\(h\in\ell_p\). Since \(\Psi(0)=0\) and \(D\Psi(0)=0\), this gives
\eqref{eq:coercive-F}.

4.)
Let \((u_n)\) be a normalized successive block sequence
and put \(k_n=\min\supp u_n\). Since
\(\norm{u_n}_\infty\le1\) and \(k_n\to\infty\),
\[
 0<Q(u_n)\le\sum_{k\ge k_n}2^{-k}=2^{1-k_n}\longrightarrow0,
 \qquad M(u_n)=1.
\]
For fixed \(t\ne0\), formula~\eqref{eq:F} therefore gives
\(\Psi(tu_n)-t^2Q(u_n)\to|t|^p\). By \eqref{eq:F-bounds},
every term of this sequence lies between zero and \(|t|^p\).
Its supremum is consequently exactly \(|t|^p\), proving
\eqref{eq:sup-remainder}. Division by \(|t|^3\), with \(p<3\),
proves the failure of a common cubic estimate.
\end{proof}

\section{Corollaries}\label{sec:cor}

\begin{corollary}\label{cor:lur}
Let \(2<p<3\). Then the following assertions hold.
\begin{enumerate}[label=\textup{\arabic*.},ref=\arabic*,leftmargin=*]
\item There is no equivalent LUR \(G^3\)-smooth norm on
\(\ell_p\).
\item More generally, no equivalent \(G^3\)-smooth norm on
\(\ell_p\) has a unit ball whose strongly exposed points are
norm-dense in its unit sphere.
\end{enumerate}
\end{corollary}
\begin{proof}
The general assertion follows from Theorem~\ref{thm:kadec}\textup{(\ref{it:kadec-exposed})}.
If \(\tnorm{\cdot}\) is LUR, a norming functional \(f\) at \(x\in S_{\tnorm{\cdot}}\) strongly
exposes \(x\): from \(f(y_n)\to1\), \(y_n\in B_{\tnorm{\cdot}}\), one obtains
\(\tnorm{y_n}\to1\); after normalization,
\(\tnorm{x+y_n/\tnorm{y_n}}\to2\), and LUR gives norm convergence.
\end{proof}

\begin{corollary}\label{cor:contact}
Let \(2<p<3\), let \(U\subset\ell_p\) be open, and let
\(f:U\to\R\) be continuous and belong to \(G^3(U)\). Suppose
there are \(a\in U\) and \(c,\delta>0\) with \(B(a,\delta)\subset U\)
such that
\begin{equation}\label{eq:contact}
 f(a+h)-f(a)-Df(a)h\ge c\norm{h}_p^p
 \qquad(\norm{h}_p<\delta).
\end{equation}
Then the following assertions hold.
\begin{enumerate}[label=\textup{\arabic*.},ref=\arabic*,leftmargin=*]
\item \(a\notin G_f(U)\).
\item All points satisfying \eqref{eq:contact} belong to
the relatively closed nowhere dense set \(U\setminus G_f(U)\).
\end{enumerate}
\end{corollary}
\begin{proof}
If \(a\in G_f(U)\), Theorem~\ref{thm:baire}\textup{(\ref{it:baire-taylor})} gives a common cubic
remainder in a ball about \(a\). Apply \eqref{eq:contact} to
\(h=t e_n\) and \(h=-t e_n\), and add the two expansions:
\[
 2c|t|^p
 \le t^2D^2f(a)[e_n,e_n]+\frac L3|t|^3.
\]
For fixed sufficiently small \(t\), let \(n\to\infty\) and use
Lemma~\ref{lem:diagonal}. Then
\(2c\le (L/3)|t|^{3-p}\), which is impossible as \(t\to0\).
\end{proof}

\begin{corollary}\label{cor:unbounded}
Let \(2<p<3\). Then the following assertions hold.
\begin{enumerate}[label=\textup{\arabic*.},ref=\arabic*,leftmargin=*]
\item If \(\tnorm{\cdot}\) is an equivalent \(G^3\)-smooth norm on
\(\ell_p\), then \(D^3\tnorm{\cdot}\) is unbounded on every
neighborhood of each strongly exposed point of \(B_{\tnorm{\cdot}}\).
\item There is no equivalent \(G^3\)-smooth norm on \(\ell_p\) whose third
derivative is locally bounded at every nonzero point.
\end{enumerate}
\end{corollary}
\begin{proof}
The first assertion is Theorem~\ref{thm:kadec}\textup{(\ref{it:kadec-regular})}.
For the second, take any strongly exposed point of the unit ball.
Such points exist because \(\ell_p\) is reflexive and its
closed bounded convex unit ball is the closed convex hull of its
strongly exposed points; see \cite[pp.~210--211]{HMZ}. Alternatively, a smooth
scalar cutoff of the norm would give a locally
\(C^{2,1}\), and hence \(T^p\)-smooth, bump, contradicting
\cite[Corollary~1.2]{DGJ}.
\end{proof}

\begin{remark}\label{rem:remaining}
Theorems~\ref{thm:baire} and \ref{thm:kadec} do not force a
strongly exposed point into \(G_{\tnorm{\cdot}}\). Reflexivity gives a closed
convex-hull representation by strongly exposed points, not their
norm-density in the sphere. Thus it does not exclude their
containment in \(S_{\tnorm{\cdot}}\setminus G_{\tnorm{\cdot}}\).
The perturbation in Lemma~\ref{lem:kadec} also explains the role of
its Kadec hypothesis. If \(x_n\rightharpoonup x\) in \(\ell_p\), then
compactness of \(T\) gives \(H(x_n)\to H(x)\). Hence
\(\tnorm{x_n}_H\to\tnorm{x}_H\) if and only if
\(\tnorm{x_n}\to\tnorm{x}\). By the sequential characterization in
Definition~\ref{def:geometry}, after normalization, the two norms
have the Kadec property simultaneously.

H\'ajek \cite{HajekKadec} proved that whenever a Banach space
admits an equivalent Kadec norm, the Kadec norms form a residual set
in the space of all equivalent norms. On \(\ell_p\), \(2<p<3\), the
standard norm is uniformly convex and hence Kadec. By
Theorem~\ref{thm:kadec}, every equivalent \(G^3\)-smooth norm belongs
to the complement of this residual set. Thus the \(G^3\)-smooth
norms form a set of first category in the space of renormings
considered in \cite{HajekKadec}. Whether this set is empty remains
open in the present argument.
\end{remark}

\section{Illustrative examples}\label{sec:examples}

\begin{example}[The canonical norm]\label{ex:canonical}
Let \(2<p<3\). Along \(e_1+t e_2\), the standard norm satisfies
\[
 \norm{e_1+t e_2}_p
 =(1+|t|^p)^{1/p}
 =1+\frac1p|t|^p+O(|t|^{2p}).
\]
Its first and second derivatives at \(t=0\) vanish, but it has no
third derivative there. Indeed, if a third derivative existed, the
evenness would make the cubic coefficient zero and the remainder
would be \(o(|t|^3)\), whereas
\[
 \frac{(1+|t|^p)^{1/p}-1}{|t|^3}
 \sim\frac1p|t|^{p-3}\longrightarrow\infty.
\]
This explains the failure for the canonical norm but does not
answer the equivalent-renorming question.
\end{example}

\begin{example}[An explicit exponent]\label{ex:five-halves}
For \(p=5/2\), choose \(m=3\) in Theorem~\ref{thm:counter}.
Then
\[
 \Psi(x)=\left[
       \left(\sum_{k=1}^{\infty}2^{-k}x_k^2\right)^3
       +\left(\sum_{k=1}^{\infty}|x_k|^{5/2}\right)^3
       \right]^{1/3}.
\]
For every fixed \(h\ne0\),
\[
 \Psi(th)=t^2Q(h)+O_h(|t|^{7/2}),
 \qquad
 \norm{D^2\Psi(th)-D^2\Psi(0)}=O_h(|t|^{3/2}).
\]
For each \(t\ne0\), every normalized successive block sequence satisfies
\[
 \sup_n\frac{|\Psi(tu_n)-t^2Q(u_n)|}{|t|^3}=|t|^{-1/2}.
\]
The constants in the fixed-direction estimates cannot therefore be
made uniform by passing to further blocks.

Coordinating directions and steps cannot make both the quadratic
term and the remainder negligible on the \(|t_n|^p\) scale.
For any \(2<p<3\), choose \(m\) as in Theorem~\ref{thm:counter}.
Given normalized blocks \(u_n\) and nonzero scalars \(t_n\to0\), put
\(a_n=Q(u_n)/|t_n|^{p-2}\). Directly from \eqref{eq:F},
\[
 \frac{\Psi(t_nu_n)-t_n^2Q(u_n)}{|t_n|^p}
 =(a_n^m+1)^{1/m}-a_n.
\]
If \(a_n\to0\), so that the quadratic term is negligible at the
scale needed in a contact argument, the remainder divided by
\(|t_n|^p\) tends to one.
\end{example}

\begin{example}[The counterexample is not globally \(G^3\)]
\label{ex:not-global}
Let \(\Psi\) be as in Theorem~\ref{thm:counter} and put
\(\Phi(a,b)=(a^m+b^m)^{1/m}\) for \(a,b\ge0\), so that
\(\Psi(x)=\Phi(Q(x),M(x))\). Since
\[
 Q(e_1+t e_2)=\tfrac12+\tfrac14t^2,
 \qquad M(e_1+t e_2)=1+|t|^p,
\]
the ordinary two-variable Taylor formula for \(\Phi\) yields
\[
 \Psi(e_1+t e_2)
 =\Psi(e_1)+\tfrac14\Phi_a(\tfrac12,1)t^2
            +\Phi_b(\tfrac12,1)|t|^p+O(|t|^4).
\]
Here \(\Phi_b(\tfrac12,1)>0\), and \(p<3\). The restriction has no
third derivative at zero. This confirms that
Theorem~\ref{thm:counter} is a counterexample to a local uniformity
claim, not a construction of a globally \(G^3\)-smooth norm or bump.
\end{example}

\begin{example}[A discontinuous second \Fre\ derivative]
\label{ex:hessian-discontinuous}
The function \(\Psi\) of Theorem~\ref{thm:counter} is twice
\Fre\ differentiable everywhere, but its second derivative
is not continuous at zero. For \(n\ge2\), put
\(t_n=2^{-n/[2(p-2)]}\) and \(x_n=t_ne_n\). Then
\(x_n\to0\), while
\(Q(x_n)/M(x_n)=2^{-n}t_n^{2-p}=2^{-n/2}\to0\).
In the direction \(e_1\), both first derivatives
\(DQ(x_n)e_1\) and \(DM(x_n)e_1\) vanish, and
\(D^2M(x_n)[e_1,e_1]=0\). Since
\(D^2Q[e_1,e_1]=1\), the second-order chain rule yields
\[
 D^2\Psi(x_n)[e_1,e_1]
 =\left(\frac{Q(x_n)}{\Psi(x_n)}\right)^{m-1}
 \longrightarrow0,
 \qquad D^2\Psi(0)[e_1,e_1]=1.
\]
Thus \(D^2\Psi(x_n)\) does not converge to \(D^2\Psi(0)\) in
operator norm. This also explains why the estimate along each
fixed direction in \eqref{eq:F-derivative-errors} must not be
read as a uniform estimate for directions varying with \(n\).
\end{example}

\begin{example}[A strictly convex smooth patch without strong exposure]
\label{ex:cone}
Fix \(2<p<3\). For \(\varepsilon>0\), define on \(\ell_p\)
\[
 \tnorm{x}_\varepsilon
 =\max\{\norm{x}_p,\,2|x_1|\}+\varepsilon\sqrt{Q(x)}.
\]
This is an equivalent strictly convex norm, by the same strict
triangle-inequality argument used in Lemma~\ref{lem:kadec}.
On the nonempty open cone
\[
 \mathcal C=\{x:2|x_1|>\norm{x}_p\}
\]
it is \(C^\infty\). Nevertheless, no point of
\(\mathcal C\cap S_{\tnorm{\cdot}_\varepsilon}\) is strongly exposed.

To see this, fix \(x\) in that set and choose a sufficiently small
fixed \(t>0\). The Lipschitz estimate
\(\norm{x+t e_n}_p\le\norm{x}_p+t\) ensures that \(x+t e_n\in
\mathcal C\) for every \(n\ge2\). Moreover,
\[
 Q(x+t e_n)-Q(x)=2^{1-n}tx_n+2^{-n}t^2\longrightarrow0.
\]
Thus \(\tnorm{x+t e_n}_\varepsilon\to1\).
Let
\[
 y_n=\frac{x+t e_n}{\tnorm{x+t e_n}_\varepsilon}.
\]
Then \(y_n\in S_{\tnorm{\cdot}_\varepsilon}\), \(y_n\rightharpoonup x\), but
\(\norm{y_n-x}_p\to t\). Every norming functional \(f\) at \(x\)
satisfies \(f(y_n)\to f(x)=1\), so it cannot strongly expose \(x\).
This example shows why strict convexity alone does not put a useful
geometric point inside a smooth open patch. No global
\(G^3\)-smoothness is asserted for \(\tnorm{\cdot}_\varepsilon\).
\end{example}

\begin{example}[The two exponent thresholds]\label{ex:thresholds}
The proof uses \(p>2\) in Lemma~\ref{lem:diagonal} and \(p<3\) in
the summability choice of Lemma~\ref{lem:bump}. At \(p=2\),
the form \(B(x,y)=\sum_kx_ky_k\) has \(B(e_n,e_n)=1\), and the
Hilbert norm is both LUR and \(C^\infty\) away from zero.
At \(p=3\), the requirements
\(\sum_jt_j^p=\infty\) and \(\sum_jt_j^3<\infty\) are incompatible.
Thus neither endpoint is covered by the argument. This is also
consistent with the distinction between \Gat\ and \Fre\
smoothness at odd integer exponents in \cite{Troyanski}.
\end{example}

The examples identify the unresolved step precisely: the global
differentiability assumption yields regular open regions, but a
general equivalent norm need not place a strongly exposed point in
those regions. The results above settle the Kadec and LUR cases
under Definition~\ref{def:gateaux}; they do not settle unrestricted
third-order \Gat\ renormability of \(\ell_p\), \(2<p<3\).

\end{document}